\documentclass[11pt]{article}

\long\def\ig#1{\relax}
\ig{Thanks to Roberto Minio for this def'n.  Compare the def'n of
\comment in AMSTeX.}

\newcount \coefa
\newcount \coefb
\newcount \coefc
\newcount\tempcounta
\newcount\tempcountb
\newcount\tempcountc
\newcount\tempcountd
\newcount\xext
\newcount\yext
\newcount\xoff
\newcount\yoff
\newcount\gap%
\newcount\arrowtypea
\newcount\arrowtypeb
\newcount\arrowtypec
\newcount\arrowtyped
\newcount\arrowtypee
\newcount\height
\newcount\width
\newcount\xpos
\newcount\ypos
\newcount\run
\newcount\rise
\newcount\arrowlength
\newcount\halflength
\newcount\arrowtype
\newdimen\tempdimen
\newdimen\xlen
\newdimen\ylen
\newsavebox{\tempboxa}%
\newsavebox{\tempboxb}%
\newsavebox{\tempboxc}%

\makeatletter
\def\settypes(#1,#2,#3){\arrowtypea#1 \arrowtypeb#2 \arrowtypec#3}
\def\settoheight#1#2{\setbox\@tempboxa\hbox{#2}#1\ht\@tempboxa\relax}%
\def\settodepth#1#2{\setbox\@tempboxa\hbox{#2}#1\dp\@tempboxa\relax}%
\def\settokens[#1`#2`#3`#4]{%
     \def\tokena{#1}\def\tokenb{#2}\def\tokenc{#3}\def\tokend{#4}}
\def\setsqparms[#1`#2`#3`#4;#5`#6]{%
\arrowtypea #1
\arrowtypeb #2
\arrowtypec #3
\arrowtyped #4
\width #5
\height #6
}
\def\setpos(#1,#2){\xpos=#1 \ypos#2}

\def\bfig{\begin{picture}(\xext,\yext)(\xoff,\yoff)}
\def\efig{\end{picture}}

\def\putbox(#1,#2)#3{\put(#1,#2){\makebox(0,0){$#3$}}}

\def\settriparms[#1`#2`#3;#4]{\settripairparms[#1`#2`#3`1`1;#4]}%

\def\settripairparms[#1`#2`#3`#4`#5;#6]{%
\arrowtypea #1
\arrowtypeb #2
\arrowtypec #3
\arrowtyped #4
\arrowtypee #5
\width #6
\height #6
}

\def\resetparms{\settripairparms[1`1`1`1`1;500]\width 500}

\resetparms

\def\mvector(#1,#2)#3{
\put(0,0){\vector(#1,#2){#3}}%
\put(0,0){\vector(#1,#2){30}}%
}
\def\evector(#1,#2)#3{{
\arrowlength #3
\put(0,0){\vector(#1,#2){\arrowlength}}%
\advance \arrowlength by-30
\put(0,0){\vector(#1,#2){\arrowlength}}%
}}

\def\horsize#1#2{%
\settowidth{\tempdimen}{$#2$}%
#1=\tempdimen
\divide #1 by\unitlength
}

\def\vertsize#1#2{%
\settoheight{\tempdimen}{$#2$}%
#1=\tempdimen
\settodepth{\tempdimen}{$#2$}%
\advance #1 by\tempdimen
\divide #1 by\unitlength
}

\def\vertadjust[#1`#2`#3]{%
\vertsize{\tempcounta}{#1}%
\vertsize{\tempcountb}{#2}%
\ifnum \tempcounta<\tempcountb \tempcounta=\tempcountb \fi
\divide\tempcounta by2
\vertsize{\tempcountb}{#3}%
\ifnum \tempcountb>0 \advance \tempcountb by20 \fi
\ifnum \tempcounta<\tempcountb \tempcounta=\tempcountb \fi
}

\def\horadjust[#1`#2`#3]{%
\horsize{\tempcounta}{#1}%
\horsize{\tempcountb}{#2}%
\ifnum \tempcounta<\tempcountb \tempcounta=\tempcountb \fi
\divide\tempcounta by20
\horsize{\tempcountb}{#3}%
\ifnum \tempcountb>0 \advance \tempcountb by60 \fi
\ifnum \tempcounta<\tempcountb \tempcounta=\tempcountb \fi
}

\ig{ In this procedure, #1 is the paramater that sticks out all the way,
#2 sticks out the least and #3 is a label sticking out half way.  #4 is
the amount of the offset.}

\def\sladjust[#1`#2`#3]#4{%
\tempcountc=#4
\horsize{\tempcounta}{#1}%
\divide \tempcounta by2
\horsize{\tempcountb}{#2}%
\divide \tempcountb by2
\advance \tempcountb by-\tempcountc
\ifnum \tempcounta<\tempcountb \tempcounta=\tempcountb\fi
\divide \tempcountc by2
\horsize{\tempcountb}{#3}%
\advance \tempcountb by-\tempcountc
\ifnum \tempcountb>0 \advance \tempcountb by80\fi
\ifnum \tempcounta<\tempcountb \tempcounta=\tempcountb\fi
\advance\tempcounta by20
}

\def\putvector(#1,#2)(#3,#4)#5#6{{%
\xpos=#1
\ypos=#2
\run=#3
\rise=#4
\arrowlength=#5
\arrowtype=#6
\ifnum \arrowtype<0
    \ifnum \run=0
        \advance \ypos by-\arrowlength
    \else
        \tempcounta \arrowlength
        \multiply \tempcounta by\rise
        \divide \tempcounta by\run
        \ifnum\run>0
            \advance \xpos by\arrowlength
            \advance \ypos by\tempcounta
        \else
            \advance \xpos by-\arrowlength
            \advance \ypos by-\tempcounta
        \fi
    \fi
    \multiply \arrowtype by-1
    \multiply \rise by-1
    \multiply \run by-1
\fi
\ifnum \arrowtype=1
    \put(\xpos,\ypos){\vector(\run,\rise){\arrowlength}}%
\else\ifnum \arrowtype=2
    \put(\xpos,\ypos){\mvector(\run,\rise)\arrowlength}%
\else\ifnum\arrowtype=3
    \put(\xpos,\ypos){\evector(\run,\rise){\arrowlength}}%
\fi\fi\fi
}}

\def\putsplitvector(#1,#2)#3#4{
\xpos #1
\ypos #2
\arrowtype #4
\halflength #3
\arrowlength #3
\gap 140
\advance \halflength by-\gap
\divide \halflength by2
\ifnum \arrowtype=1
    \put(\xpos,\ypos){\line(0,-1){\halflength}}%
    \advance\ypos by-\halflength
    \advance\ypos by-\gap
    \put(\xpos,\ypos){\vector(0,-1){\halflength}}%
\else\ifnum \arrowtype=2
    \put(\xpos,\ypos){\line(0,-1)\halflength}%
    \put(\xpos,\ypos){\vector(0,-1)3}%
    \advance\ypos by-\halflength
    \advance\ypos by-\gap
    \put(\xpos,\ypos){\vector(0,-1){\halflength}}%
\else\ifnum\arrowtype=3
    \put(\xpos,\ypos){\line(0,-1)\halflength}%
    \advance\ypos by-\halflength
    \advance\ypos by-\gap
    \put(\xpos,\ypos){\evector(0,-1){\halflength}}%
\else\ifnum \arrowtype=-1
    \advance \ypos by-\arrowlength
    \put(\xpos,\ypos){\line(0,1){\halflength}}%
    \advance\ypos by\halflength
    \advance\ypos by\gap
    \put(\xpos,\ypos){\vector(0,1){\halflength}}%
\else\ifnum \arrowtype=-2
    \advance \ypos by-\arrowlength
    \put(\xpos,\ypos){\line(0,1)\halflength}%
    \put(\xpos,\ypos){\vector(0,1)3}%
    \advance\ypos by\halflength
    \advance\ypos by\gap
    \put(\xpos,\ypos){\vector(0,1){\halflength}}%
\else\ifnum\arrowtype=-3
    \advance \ypos by-\arrowlength
    \put(\xpos,\ypos){\line(0,1)\halflength}%
    \advance\ypos by\halflength
    \advance\ypos by\gap
    \put(\xpos,\ypos){\evector(0,1){\halflength}}%
\fi\fi\fi\fi\fi\fi
}

\def\putmorphism(#1)(#2,#3)[#4`#5`#6]#7#8#9{{%
\run #2
\rise #3
\ifnum\rise=0
  \puthmorphism(#1)[#4`#5`#6]{#7}{#8}{#9}%
\else\ifnum\run=0
  \putvmorphism(#1)[#4`#5`#6]{#7}{#8}{#9}%
\else
\setpos(#1)%
\arrowlength #7
\arrowtype #8
\ifnum\run=0
\else\ifnum\rise=0
\else
\ifnum\run>0
    \coefa=1
\else
   \coefa=-1
\fi
\ifnum\arrowtype>0
   \coefb=0
   \coefc=-1
\else
   \coefb=\coefa
   \coefc=1
   \arrowtype=-\arrowtype
\fi
\width=2
\multiply \width by\run
\divide \width by\rise
\ifnum \width<0  \width=-\width\fi
\advance\width by60
\if l#9 \width=-\width\fi
\putbox(\xpos,\ypos){#4}
{\multiply \coefa by\arrowlength
\advance\xpos by\coefa
\multiply \coefa by\rise
\divide \coefa by\run
\advance \ypos by\coefa
\putbox(\xpos,\ypos){#5} }%
{\multiply \coefa by\arrowlength
\divide \coefa by2
\advance \xpos by\coefa
\advance \xpos by\width
\multiply \coefa by\rise
\divide \coefa by\run
\advance \ypos by\coefa
\if l#9%
   \put(\xpos,\ypos){\makebox(0,0)[r]{$#6$}}%
\else\if r#9%
   \put(\xpos,\ypos){\makebox(0,0)[l]{$#6$}}%
\fi\fi }%
{\multiply \rise by-\coefc
\multiply \run by-\coefc
\multiply \coefb by\arrowlength
\advance \xpos by\coefb
\multiply \coefb by\rise
\divide \coefb by\run
\advance \ypos by\coefb
\multiply \coefc by70
\advance \ypos by\coefc
\multiply \coefc by\run
\divide \coefc by\rise
\advance \xpos by\coefc
\multiply \coefa by140
\multiply \coefa by\run
\divide \coefa by\rise
\advance \arrowlength by\coefa
\ifnum \arrowtype=1
   \put(\xpos,\ypos){\vector(\run,\rise){\arrowlength}}%
\else\ifnum\arrowtype=2
   \put(\xpos,\ypos){\mvector(\run,\rise){\arrowlength}}%
\else\ifnum\arrowtype=3
   \put(\xpos,\ypos){\evector(\run,\rise){\arrowlength}}%
\fi\fi\fi}\fi\fi\fi\fi}}

\def\puthmorphism(#1,#2)[#3`#4`#5]#6#7#8{{%
\xpos #1
\ypos #2
\width #6
\arrowlength #6
\putbox(\xpos,\ypos){#3\vphantom{#4}}%
{\advance \xpos by\arrowlength
\putbox(\xpos,\ypos){\vphantom{#3}#4}}%
\horsize{\tempcounta}{#3}%
\horsize{\tempcountb}{#4}%
\divide \tempcounta by2
\divide \tempcountb by2
\advance \tempcounta by30
\advance \tempcountb by30
\advance \xpos by\tempcounta
\advance \arrowlength by-\tempcounta
\advance \arrowlength by-\tempcountb
\putvector(\xpos,\ypos)(1,0){\arrowlength}{#7}%
\divide \arrowlength by2
\advance \xpos by\arrowlength
\vertsize{\tempcounta}{#5}%
\divide\tempcounta by2
\advance \tempcounta by20
\if a#8 %
   \advance \ypos by\tempcounta
   \putbox(\xpos,\ypos){#5}%
\else
   \advance \ypos by-\tempcounta
   \putbox(\xpos,\ypos){#5}%
\fi}}

\def\putvmorphism(#1,#2)[#3`#4`#5]#6#7#8{{%
\xpos #1
\ypos #2
\arrowlength #6
\arrowtype #7
\settowidth{\xlen}{$#5$}%
\putbox(\xpos,\ypos){#3}%
{\advance \ypos by-\arrowlength
\putbox(\xpos,\ypos){#4}}%
{\advance\arrowlength by-140
\advance \ypos by-70
\ifdim\xlen>0pt
   \if m#8%
      \putsplitvector(\xpos,\ypos){\arrowlength}{\arrowtype}%
   \else
      \putvector(\xpos,\ypos)(0,-1){\arrowlength}{\arrowtype}%
   \fi
\else
   \putvector(\xpos,\ypos)(0,-1){\arrowlength}{\arrowtype}%
\fi}%
\ifdim\xlen>0pt
   \divide \arrowlength by2
   \advance\ypos by-\arrowlength
   \if l#8%
      \advance \xpos by-40
      \put(\xpos,\ypos){\makebox(0,0)[r]{$#5$}}%
   \else\if r#8%
      \advance \xpos by40
      \put(\xpos,\ypos){\makebox(0,0)[l]{$#5$}}%
   \else
      \putbox(\xpos,\ypos){#5}%
   \fi\fi
\fi
}}

\def\topadjust[#1`#2`#3]{%
\yoff=10
\vertadjust[#1`#2`{#3}]%
\advance \yext by\tempcounta
\advance \yext by 10
}
\def\botadjust[#1`#2`#3]{%
\vertadjust[#1`#2`{#3}]%
\advance \yext by\tempcounta
\advance \yoff by-\tempcounta
}
\def\leftadjust[#1`#2`#3]{%
\xoff=0
\horadjust[#1`#2`{#3}]%
\advance \xext by\tempcounta
\advance \xoff by-\tempcounta
}
\def\rightadjust[#1`#2`#3]{%
\horadjust[#1`#2`{#3}]%
\advance \xext by\tempcounta
}
\def\rightsladjust[#1`#2`#3]{%
\sladjust[#1`#2`{#3}]{\width}%
\advance \xext by\tempcounta
}
\def\leftsladjust[#1`#2`#3]{%
\xoff=0
\sladjust[#1`#2`{#3}]{\width}%
\advance \xext by\tempcounta
\advance \xoff by-\tempcounta
}
\def\adjust[#1`#2;#3`#4;#5`#6;#7`#8]{%
\topadjust[#1``{#2}]
\leftadjust[#3``{#4}]
\rightadjust[#5``{#6}]
\botadjust[#7``{#8}]}

\def\putsquarep<#1>(#2)[#3;#4`#5`#6`#7]{{%
\setsqparms[#1]%
\setpos(#2)%
\settokens[#3]%
\puthmorphism(\xpos,\ypos)[\tokenc`\tokend`{#7}]{\width}{\arrowtyped}b%
\advance\ypos by \height
\puthmorphism(\xpos,\ypos)[\tokena`\tokenb`{#4}]{\width}{\arrowtypea}a%
\putvmorphism(\xpos,\ypos)[``{#5}]{\height}{\arrowtypeb}l%
\advance\xpos by \width
\putvmorphism(\xpos,\ypos)[``{#6}]{\height}{\arrowtypec}r%
}}

\def\putsquare{\@ifnextchar <{\putsquarep}{\putsquarep%
   <\arrowtypea`\arrowtypeb`\arrowtypec`\arrowtyped;\width`\height>}}
\def\square{\@ifnextchar< {\squarep}{\squarep
   <\arrowtypea`\arrowtypeb`\arrowtypec`\arrowtyped;\width`\height>}}
\def\squarep<#1>[#2`#3`#4`#5;#6`#7`#8`#9]{{
\setsqparms[#1]
\xext=\width                                          
\yext=\height                                         
\topadjust[#2`#3`{#6}]
\botadjust[#4`#5`{#9}]
\leftadjust[#2`#4`{#7}]
\rightadjust[#3`#5`{#8}]
\begin{picture}(\xext,\yext)(\xoff,\yoff)
\putsquarep<\arrowtypea`\arrowtypeb`\arrowtypec`\arrowtyped;\width`\height>%
(0,0)[#2`#3`#4`#5;#6`#7`#8`{#9}]%
\end{picture}%
}}

\def\putptrianglep<#1>(#2,#3)[#4`#5`#6;#7`#8`#9]{{%
\settriparms[#1]%
\xpos=#2 \ypos=#3
\advance\ypos by \height
\puthmorphism(\xpos,\ypos)[#4`#5`{#7}]{\height}{\arrowtypea}a%
\putvmorphism(\xpos,\ypos)[`#6`{#8}]{\height}{\arrowtypeb}l%
\advance\xpos by\height
\putmorphism(\xpos,\ypos)(-1,-1)[``{#9}]{\height}{\arrowtypec}r%
}}

\def\putptriangle{\@ifnextchar <{\putptrianglep}{\putptrianglep
   <\arrowtypea`\arrowtypeb`\arrowtypec;\height>}}
\def\ptriangle{\@ifnextchar <{\ptrianglep}{\ptrianglep
   <\arrowtypea`\arrowtypeb`\arrowtypec;\height>}}

\def\ptrianglep<#1>[#2`#3`#4;#5`#6`#7]{{
\settriparms[#1]%
\width=\height                         
\xext=\width                           
\yext=\width                           
\topadjust[#2`#3`{#5}]
\botadjust[#3``]
\leftadjust[#2`#4`{#6}]
\rightsladjust[#3`#4`{#7}]
\begin{picture}(\xext,\yext)(\xoff,\yoff)
\putptrianglep<\arrowtypea`\arrowtypeb`\arrowtypec;\height>%
(0,0)[#2`#3`#4;#5`#6`{#7}]%
\end{picture}%
}}

\def\putqtrianglep<#1>(#2,#3)[#4`#5`#6;#7`#8`#9]{{%
\settriparms[#1]%
\xpos=#2 \ypos=#3
\advance\ypos by\height
\puthmorphism(\xpos,\ypos)[#4`#5`{#7}]{\height}{\arrowtypea}a%
\putmorphism(\xpos,\ypos)(1,-1)[``{#8}]{\height}{\arrowtypeb}l%
\advance\xpos by\height
\putvmorphism(\xpos,\ypos)[`#6`{#9}]{\height}{\arrowtypec}r%
}}

\def\putqtriangle{\@ifnextchar <{\putqtrianglep}{\putqtrianglep
   <\arrowtypea`\arrowtypeb`\arrowtypec;\height>}}
\def\qtriangle{\@ifnextchar <{\qtrianglep}{\qtrianglep
   <\arrowtypea`\arrowtypeb`\arrowtypec;\height>}}

\def\qtrianglep<#1>[#2`#3`#4;#5`#6`#7]{{
\settriparms[#1]
\width=\height                         
\xext=\width                           
\yext=\height                          
\topadjust[#2`#3`{#5}]
\botadjust[#4``]
\leftsladjust[#2`#4`{#6}]
\rightadjust[#3`#4`{#7}]
\begin{picture}(\xext,\yext)(\xoff,\yoff)
\putqtrianglep<\arrowtypea`\arrowtypeb`\arrowtypec;\height>%
(0,0)[#2`#3`#4;#5`#6`{#7}]%
\end{picture}%
}}

\def\putdtrianglep<#1>(#2,#3)[#4`#5`#6;#7`#8`#9]{{%
\settriparms[#1]%
\xpos=#2 \ypos=#3
\puthmorphism(\xpos,\ypos)[#5`#6`{#9}]{\height}{\arrowtypec}b%
\advance\xpos by \height \advance\ypos by\height
\putmorphism(\xpos,\ypos)(-1,-1)[``{#7}]{\height}{\arrowtypea}l%
\putvmorphism(\xpos,\ypos)[#4``{#8}]{\height}{\arrowtypeb}r%
}}

\def\putdtriangle{\@ifnextchar <{\putdtrianglep}{\putdtrianglep
   <\arrowtypea`\arrowtypeb`\arrowtypec;\height>}}
\def\dtriangle{\@ifnextchar <{\dtrianglep}{\dtrianglep
   <\arrowtypea`\arrowtypeb`\arrowtypec;\height>}}

\def\dtrianglep<#1>[#2`#3`#4;#5`#6`#7]{{
\settriparms[#1]
\width=\height                         
\xext=\width                           
\yext=\height                          
\topadjust[#2``]
\botadjust[#3`#4`{#7}]
\leftsladjust[#3`#2`{#5}]
\rightadjust[#2`#4`{#6}]
\begin{picture}(\xext,\yext)(\xoff,\yoff)
\putdtrianglep<\arrowtypea`\arrowtypeb`\arrowtypec;\height>%
(0,0)[#2`#3`#4;#5`#6`{#7}]%
\end{picture}%
}}

\def\putbtrianglep<#1>(#2,#3)[#4`#5`#6;#7`#8`#9]{{%
\settriparms[#1]%
\xpos=#2 \ypos=#3
\puthmorphism(\xpos,\ypos)[#5`#6`{#9}]{\height}{\arrowtypec}b%
\advance\ypos by\height
\putmorphism(\xpos,\ypos)(1,-1)[``{#8}]{\height}{\arrowtypeb}r%
\putvmorphism(\xpos,\ypos)[#4``{#7}]{\height}{\arrowtypea}l%
}}

\def\putbtriangle{\@ifnextchar <{\putbtrianglep}{\putbtrianglep
   <\arrowtypea`\arrowtypeb`\arrowtypec;\height>}}
\def\btriangle{\@ifnextchar <{\btrianglep}{\btrianglep
   <\arrowtypea`\arrowtypeb`\arrowtypec;\height>}}

\def\btrianglep<#1>[#2`#3`#4;#5`#6`#7]{{
\settriparms[#1]
\width=\height                         
\xext=\width                           
\yext=\height                          
\topadjust[#2``]
\botadjust[#3`#4`{#7}]
\leftadjust[#2`#3`{#5}]
\rightsladjust[#4`#2`{#6}]
\begin{picture}(\xext,\yext)(\xoff,\yoff)
\putbtrianglep<\arrowtypea`\arrowtypeb`\arrowtypec;\height>%
(0,0)[#2`#3`#4;#5`#6`{#7}]%
\end{picture}%
}}

\def\putAtrianglep<#1>(#2,#3)[#4`#5`#6;#7`#8`#9]{{%
\settriparms[#1]%
\xpos=#2 \ypos=#3
{\multiply \height by2
\puthmorphism(\xpos,\ypos)[#5`#6`{#9}]{\height}{\arrowtypec}b}%
\advance\xpos by\height \advance\ypos by\height
\putmorphism(\xpos,\ypos)(-1,-1)[#4``{#7}]{\height}{\arrowtypea}l%
\putmorphism(\xpos,\ypos)(1,-1)[``{#8}]{\height}{\arrowtypeb}r%
}}

\def\putAtriangle{\@ifnextchar <{\putAtrianglep}{\putAtrianglep
   <\arrowtypea`\arrowtypeb`\arrowtypec;\height>}}
\def\Atriangle{\@ifnextchar <{\Atrianglep}{\Atrianglep
   <\arrowtypea`\arrowtypeb`\arrowtypec;\height>}}

\def\Atrianglep<#1>[#2`#3`#4;#5`#6`#7]{{
\settriparms[#1]
\width=\height                         
\xext=\width                           
\yext=\height                          
\topadjust[#2``]
\botadjust[#3`#4`{#7}]
\multiply \xext by2 
\leftsladjust[#3`#2`{#5}]
\rightsladjust[#4`#2`{#6}]
\begin{picture}(\xext,\yext)(\xoff,\yoff)%
\putAtrianglep<\arrowtypea`\arrowtypeb`\arrowtypec;\height>%
(0,0)[#2`#3`#4;#5`#6`{#7}]%
\end{picture}%
}}

\def\putAtrianglepairp<#1>(#2)[#3;#4`#5`#6`#7`#8]{{
\settripairparms[#1]%
\setpos(#2)%
\settokens[#3]%
\puthmorphism(\xpos,\ypos)[\tokenb`\tokenc`{#7}]{\height}{\arrowtyped}b%
\advance\xpos by\height
\advance\ypos by\height
\putmorphism(\xpos,\ypos)(-1,-1)[\tokena``{#4}]{\height}{\arrowtypea}l%
\putvmorphism(\xpos,\ypos)[``{#5}]{\height}{\arrowtypeb}m%
\putmorphism(\xpos,\ypos)(1,-1)[``{#6}]{\height}{\arrowtypec}r%
}}

\def\putAtrianglepair{\@ifnextchar <{\putAtrianglepairp}{\putAtrianglepairp%
   <\arrowtypea`\arrowtypeb`\arrowtypec`\arrowtyped`\arrowtypee;\height>}}
\def\Atrianglepair{\@ifnextchar <{\Atrianglepairp}{\Atrianglepairp%
   <\arrowtypea`\arrowtypeb`\arrowtypec`\arrowtyped`\arrowtypee;\height>}}

\def\Atrianglepairp<#1>[#2;#3`#4`#5`#6`#7]{{%
\settripairparms[#1]%
\settokens[#2]%
\width=\height
\xext=\width
\yext=\height
\topadjust[\tokena``]%
\vertadjust[\tokenb`\tokenc`{#6}]
\tempcountd=\tempcounta                       
\vertadjust[\tokenc`\tokend`{#7}]
\ifnum\tempcounta<\tempcountd                 
\tempcounta=\tempcountd\fi                    
\advance \yext by\tempcounta                  
\advance \yoff by-\tempcounta                 %
\multiply \xext by2 
\leftsladjust[\tokenb`\tokena`{#3}]
\rightsladjust[\tokend`\tokena`{#5}]%
\begin{picture}(\xext,\yext)(\xoff,\yoff)%
\putAtrianglepairp
<\arrowtypea`\arrowtypeb`\arrowtypec`\arrowtyped`\arrowtypee;\height>%
(0,0)[#2;#3`#4`#5`#6`{#7}]%
\end{picture}%
}}

\def\putVtrianglep<#1>(#2,#3)[#4`#5`#6;#7`#8`#9]{{%
\settriparms[#1]%
\xpos=#2 \ypos=#3
\advance\ypos by\height
{\multiply\height by2
\puthmorphism(\xpos,\ypos)[#4`#5`{#7}]{\height}{\arrowtypea}a}%
\putmorphism(\xpos,\ypos)(1,-1)[`#6`{#8}]{\height}{\arrowtypeb}l%
\advance\xpos by\height
\advance\xpos by\height
\putmorphism(\xpos,\ypos)(-1,-1)[``{#9}]{\height}{\arrowtypec}r%
}}

\def\putVtriangle{\@ifnextchar <{\putVtrianglep}{\putVtrianglep
   <\arrowtypea`\arrowtypeb`\arrowtypec;\height>}}
\def\Vtriangle{\@ifnextchar <{\Vtrianglep}{\Vtrianglep
   <\arrowtypea`\arrowtypeb`\arrowtypec;\height>}}

\def\Vtrianglep<#1>[#2`#3`#4;#5`#6`#7]{{
\settriparms[#1]
\width=\height                         
\xext=\width                           
\yext=\height                          
\topadjust[#2`#3`{#5}]
\botadjust[#4``]
\multiply \xext by2 
\leftsladjust[#2`#3`{#6}]
\rightsladjust[#3`#4`{#7}]
\begin{picture}(\xext,\yext)(\xoff,\yoff)%
\putVtrianglep<\arrowtypea`\arrowtypeb`\arrowtypec;\height>%
(0,0)[#2`#3`#4;#5`#6`{#7}]%
\end{picture}%
}}

\def\putVtrianglepairp<#1>(#2)[#3;#4`#5`#6`#7`#8]{{
\settripairparms[#1]%
\setpos(#2)%
\settokens[#3]%
\advance\ypos by\height
\putmorphism(\xpos,\ypos)(1,-1)[`\tokend`{#6}]{\height}{\arrowtypec}l%
\puthmorphism(\xpos,\ypos)[\tokena`\tokenb`{#4}]{\height}{\arrowtypea}a%
\advance\xpos by\height
\putvmorphism(\xpos,\ypos)[``{#7}]{\height}{\arrowtyped}m%
\advance\xpos by\height
\putmorphism(\xpos,\ypos)(-1,-1)[``{#8}]{\height}{\arrowtypee}r%
}}

\def\putVtrianglepair{\@ifnextchar <{\putVtrianglepairp}{\putVtrianglepairp%
    <\arrowtypea`\arrowtypeb`\arrowtypec`\arrowtyped`\arrowtypee;\height>}}
\def\Vtrianglepair{\@ifnextchar <{\Vtrianglepairp}{\Vtrianglepairp%
    <\arrowtypea`\arrowtypeb`\arrowtypec`\arrowtyped`\arrowtypee;\height>}}

\def\Vtrianglepairp<#1>[#2;#3`#4`#5`#6`#7]{{%
\settripairparms[#1]%
\settokens[#2]
\xext=\height                  
\width=\height                 
\yext=\height                  
\vertadjust[\tokena`\tokenb`{#4}]
\tempcountd=\tempcounta        
\vertadjust[\tokenb`\tokenc`{#5}]
\ifnum\tempcounta<\tempcountd%
\tempcounta=\tempcountd\fi
\advance \yext by\tempcounta
\botadjust[\tokend``]%
\multiply \xext by2
\leftsladjust[\tokena`\tokend`{#6}]%
\rightsladjust[\tokenc`\tokend`{#7}]%
\begin{picture}(\xext,\yext)(\xoff,\yoff)%
\putVtrianglepairp
<\arrowtypea`\arrowtypeb`\arrowtypec`\arrowtyped`\arrowtypee;\height>%
(0,0)[#2;#3`#4`#5`#6`{#7}]%
\end{picture}%
}}

\def\putCtrianglep<#1>(#2,#3)[#4`#5`#6;#7`#8`#9]{{%
\settriparms[#1]%
\xpos=#2 \ypos=#3
\advance\ypos by\height
\putmorphism(\xpos,\ypos)(1,-1)[``{#9}]{\height}{\arrowtypec}l%
\advance\xpos by\height
\advance\ypos by\height
\putmorphism(\xpos,\ypos)(-1,-1)[#4`#5`{#7}]{\height}{\arrowtypea}l%
{\multiply\height by 2
\putvmorphism(\xpos,\ypos)[`#6`{#8}]{\height}{\arrowtypeb}r}%
}}

\def\putCtriangle{\@ifnextchar <{\putCtrianglep}{\putCtrianglep
    <\arrowtypea`\arrowtypeb`\arrowtypec;\height>}}
\def\Ctriangle{\@ifnextchar <{\Ctrianglep}{\Ctrianglep
    <\arrowtypea`\arrowtypeb`\arrowtypec;\height>}}

\def\Ctrianglep<#1>[#2`#3`#4;#5`#6`#7]{{
\settriparms[#1]
\width=\height                          
\xext=\width                            
\yext=\height                           
\multiply \yext by2 
\topadjust[#2``]
\botadjust[#4``]
\sladjust[#3`#2`{#5}]{\width}
\tempcountd=\tempcounta                 
\sladjust[#3`#4`{#7}]{\width}
\ifnum \tempcounta<\tempcountd          
\tempcounta=\tempcountd\fi              
\advance \xext by\tempcounta            
\advance \xoff by-\tempcounta           %
\rightadjust[#2`#4`{#6}]
\begin{picture}(\xext,\yext)(\xoff,\yoff)%
\putCtrianglep<\arrowtypea`\arrowtypeb`\arrowtypec;\height>%
(0,0)[#2`#3`#4;#5`#6`{#7}]%
\end{picture}%
}}

\def\putDtrianglep<#1>(#2,#3)[#4`#5`#6;#7`#8`#9]{{%
\settriparms[#1]%
\xpos=#2 \ypos=#3
\advance\xpos by\height \advance\ypos by\height
\putmorphism(\xpos,\ypos)(-1,-1)[``{#9}]{\height}{\arrowtypec}r%
\advance\xpos by-\height \advance\ypos by\height
\putmorphism(\xpos,\ypos)(1,-1)[`#5`{#8}]{\height}{\arrowtypeb}r%
{\multiply\height by 2
\putvmorphism(\xpos,\ypos)[#4`#6`{#7}]{\height}{\arrowtypea}l}%
}}

\def\putDtriangle{\@ifnextchar <{\putDtrianglep}{\putDtrianglep
    <\arrowtypea`\arrowtypeb`\arrowtypec;\height>}}
\def\Dtriangle{\@ifnextchar <{\Dtrianglep}{\Dtrianglep
   <\arrowtypea`\arrowtypeb`\arrowtypec;\height>}}

\def\Dtrianglep<#1>[#2`#3`#4;#5`#6`#7]{{
\settriparms[#1]
\width=\height                         
\xext=\height                          
\yext=\height                          
\multiply \yext by2 
\topadjust[#2``]
\botadjust[#4``]
\leftadjust[#2`#4`{#5}]
\sladjust[#3`#2`{#5}]{\height}
\tempcountd=\tempcountd                
\sladjust[#3`#4`{#7}]{\height}
\ifnum \tempcounta<\tempcountd         
\tempcounta=\tempcountd\fi             
\advance \xext by\tempcounta           %
\begin{picture}(\xext,\yext)(\xoff,\yoff)
\putDtrianglep<\arrowtypea`\arrowtypeb`\arrowtypec;\height>%
(0,0)[#2`#3`#4;#5`#6`{#7}]%
\end{picture}%
}}

\def\setrecparms[#1`#2]{\width=#1 \height=#2}%
\def\recursep<#1`#2>[#3;#4`#5`#6`#7`#8]{{%
\width=#1 \height=#2
\settokens[#3]
\settowidth{\tempdimen}{$\tokena$}
\ifdim\tempdimen=0pt
  \savebox{\tempboxa}{\hbox{$\tokenb$}}%
  \savebox{\tempboxb}{\hbox{$\tokend$}}%
  \savebox{\tempboxc}{\hbox{$#6$}}%
\else
  \savebox{\tempboxa}{\hbox{$\hbox{$\tokena$}\times\hbox{$\tokenb$}$}}%
  \savebox{\tempboxb}{\hbox{$\hbox{$\tokena$}\times\hbox{$\tokend$}$}}%
  \savebox{\tempboxc}{\hbox{$\hbox{$\tokena$}\times\hbox{$#6$}$}}%
\fi
\ypos=\height
\divide\ypos by 2
\xpos=\ypos
\advance\xpos by \width
\xext=\xpos \yext=\height
\topadjust[#3`\usebox{\tempboxa}`{#4}]%
\botadjust[#5`\usebox{\tempboxb}`{#8}]%
\sladjust[\tokenc`\tokenb`{#5}]{\ypos}%
\tempcountd=\tempcounta
\sladjust[\tokenc`\tokend`{#5}]{\ypos}%
\ifnum \tempcounta<\tempcountd
\tempcounta=\tempcountd\fi
\advance \xext by\tempcounta
\advance \xoff by-\tempcounta
\rightadjust[\usebox{\tempboxa}`\usebox{\tempboxb}`\usebox{\tempboxc}]%
\bfig
\putCtrianglep<-1`1`1;\ypos>(0,0)[`\tokenc`;#5`#6`{#7}]%
\puthmorphism(\ypos,0)[\tokend`\usebox{\tempboxb}`{#8}]{\width}{-1}b%
\puthmorphism(\ypos,\height)[\tokenb`\usebox{\tempboxa}`{#4}]{\width}{-1}a%
\advance\ypos by \width
\putvmorphism(\ypos,\height)[``\usebox{\tempboxc}]{\height}1r%
\efig
}}

\def\recurse{\@ifnextchar <{\recursep}{\recursep<\width`\height>}}

\def\puttwohmorphisms(#1,#2)[#3`#4;#5`#6]#7#8#9{{%
\puthmorphism(#1,#2)[#3`#4`]{#7}0a
\ypos=#2
\advance\ypos by 20
\puthmorphism(#1,\ypos)[\phantom{#3}`\phantom{#4}`#5]{#7}{#8}a
\advance\ypos by -40
\puthmorphism(#1,\ypos)[\phantom{#3}`\phantom{#4}`#6]{#7}{#9}b
}}

\def\puttwovmorphisms(#1,#2)[#3`#4;#5`#6]#7#8#9{{%
\putvmorphism(#1,#2)[#3`#4`]{#7}0a
\xpos=#1
\advance\xpos by -20
\putvmorphism(\xpos,#2)[\phantom{#3}`\phantom{#4}`#5]{#7}{#8}l
\advance\xpos by 40
\putvmorphism(\xpos,#2)[\phantom{#3}`\phantom{#4}`#6]{#7}{#9}r
}}

\def\puthcoequalizer(#1)[#2`#3`#4;#5`#6`#7]#8#9{{%
\setpos(#1)%
\puttwohmorphisms(\xpos,\ypos)[#2`#3;#5`#6]{#8}11%
\advance\xpos by #8
\puthmorphism(\xpos,\ypos)[\phantom{#3}`#4`#7]{#8}1{#9}
}}

\def\putvcoequalizer(#1)[#2`#3`#4;#5`#6`#7]#8#9{{%
\setpos(#1)%
\puttwovmorphisms(\xpos,\ypos)[#2`#3;#5`#6]{#8}11%
\advance\ypos by -#8
\putvmorphism(\xpos,\ypos)[\phantom{#3}`#4`#7]{#8}1{#9}
}}

\def\putthreehmorphisms(#1)[#2`#3;#4`#5`#6]#7(#8)#9{{%
\setpos(#1) \settypes(#8)
\if a#9 %
     \vertsize{\tempcounta}{#5}%
     \vertsize{\tempcountb}{#6}%
     \ifnum \tempcounta<\tempcountb \tempcounta=\tempcountb \fi
\else
     \vertsize{\tempcounta}{#4}%
     \vertsize{\tempcountb}{#5}%
     \ifnum \tempcounta<\tempcountb \tempcounta=\tempcountb \fi
\fi
\advance \tempcounta by 60
\puthmorphism(\xpos,\ypos)[#2`#3`#5]{#7}{\arrowtypeb}{#9}
\advance\ypos by \tempcounta
\puthmorphism(\xpos,\ypos)[\phantom{#2}`\phantom{#3}`#4]{#7}{\arrowtypea}{#9}
\advance\ypos by -\tempcounta \advance\ypos by -\tempcounta
\puthmorphism(\xpos,\ypos)[\phantom{#2}`\phantom{#3}`#6]{#7}{\arrowtypec}{#9}
}}

\def\putarc(#1,#2)[#3`#4`#5]#6#7#8{{%
\xpos #1
\ypos #2
\width #6
\arrowlength #6
\putbox(\xpos,\ypos){#3\vphantom{#4}}%
{\advance \xpos by\arrowlength
\putbox(\xpos,\ypos){\vphantom{#3}#4}}%
\horsize{\tempcounta}{#3}%
\horsize{\tempcountb}{#4}%
\divide \tempcounta by2
\divide \tempcountb by2
\advance \tempcounta by30
\advance \tempcountb by30
\advance \xpos by\tempcounta
\advance \arrowlength by-\tempcounta
\advance \arrowlength by-\tempcountb
\halflength=\arrowlength \divide\halflength by 2
\divide\arrowlength by 5
\put(\xpos,\ypos){\bezier{\arrowlength}(0,0)(50,50)(\halflength,50)}
\ifnum #7=-1 \put(\xpos,\ypos){\vector(-3,-2)0} \fi
\advance\xpos by \halflength
\put(\xpos,\ypos){\xpos=\halflength \advance\xpos by -50
   \bezier{\arrowlength}(0,50)(\xpos,50)(\halflength,0)}
\ifnum #7=1 {\advance \xpos by
   \halflength \put(\xpos,\ypos){\vector(3,-2)0}} \fi
\advance\ypos by 50
\vertsize{\tempcounta}{#5}%
\divide\tempcounta by2
\advance \tempcounta by20
\if a#8 %
   \advance \ypos by\tempcounta
   \putbox(\xpos,\ypos){#5}%
\else
   \advance \ypos by-\tempcounta
   \putbox(\xpos,\ypos){#5}%
\fi
}}

\makeatother

\usepackage{dsfont}
\usepackage{stmaryrd}
\usepackage{enumerate}
\usepackage{hyperref}

\def\env{{\bf env} }
\def\bv{{\bf bv} }
\def\var{{\bf var} }

\def\lk{\langle}
\def\rk{\rangle}

\newcommand{\lra}{\longrightarrow}

\newcommand{\ra}{\rightarrow}

\def\cP{{\cal P}}

\def\T{\mathds{T}}

\begin{document}

\title{Generalized Quantifiers on Dependent Types:\\
A System for Anaphora}
\author{Justyna Grudzi\'{n}ska and Marek Zawadowski}

\date{}

\maketitle

\noindent \textbf{Abstract} We propose a system for the interpretation of anaphoric relationships between unbound pronouns and quantifiers. The main technical contribution of our proposal consists in combining generalized quantifiers (\cite{MA}, \cite{LP}, \cite{BC}) with dependent types (\cite{MLP}, \cite{RA}, \cite{MM}). Empirically, our system allows a uniform treatment of the major types of unbound anaphora, with the anaphoric (dynamic) effects falling out naturally as a consequence of having generalized quantification on dependent types.


\section{Unbound anaphora} \label{4}

   A fundamental insight of dynamic semantics is that quantificational sentences have the ability to change contexts by setting up new referents (e.g., sets, dependencies) and anaphoric pronouns have the ability to refer back to them (\cite{KR}, \cite{vanB}). This paper proposes a uniform mechanism to account for a wide range of anaphoric (dynamic) effects associated with natural language quantification

\begin{itemize}
  \item Maximal anaphora to quantifiers \\
  E.g.: \underline{Most kids} entered. \underline{They} looked happy.\\
  The observation in \cite{KR}, \cite{vanB}, \cite{NR} is that the anaphoric pronoun \textit{they} in the second sentence (what we will call an anaphoric continuation) refers to the entire set of kids who entered. Thus the first sentence must introduce the set of all kids who entered.
  \item Quantificational subordination\\
  E.g.: \underline{Every man} loves \underline{a woman}. \underline{They} (each) kiss \underline{them}.\\
  The most obvious way to understand the anaphoric continuation is that every man kisses the women he loves rather than those loved by someone else (\cite{KR}, \cite{KM}, \cite{vanB}, \cite{NR}). Thus the first sentence must introduce a dependency between each of the men and the women they love that can be elaborated upon in further discourse
  \item Cumulative and branching continuations \\
  E.g.: Last year \underline{three scientists} wrote \underline{five papers}. \underline{They} presented \underline{them} at major conferences.\\
  The first sentence allows the so-called cumulative and branching readings. On the cumulative reading, it is understood to mean: Last year three scientists wrote (a total of) five papers (between them). On the branching reading, it is understood to mean: Last year three scientists (each) wrote (the same) five papers. The observation in \cite{KM}, \cite{DP08} is that the dynamics of the first sentence can deliver some cumulative or branching relation that can be elaborated upon in the anaphoric continuation .
  \item `Donkey anaphora'\\
  E.g.: Every \underline{farmer who owns a donkey} beats \underline{it}.\\
  This example shows that context can get changed within a single sentence itself (\cite{KH}, \cite{GS}, \cite{KR}). Here the modified common nouns (e.g., \textit{farmer who owns a donkey}) must introduce referents (possibly dependencies) for the respective pronouns to pick up.
\end{itemize}
The phenomenon is known as `unbound anaphora', as it refers to instances where anaphoric pronouns occur outside the syntactic scopes (i.e. the c-command domain) of their quantifier antecedents - the anaphoric pronouns are not syntactically bound by their quantifier antecedents. Unbound anaphora has been dealt with in three main semantic paradigms
\begin{itemize}
  \item Dynamic semantic theories (\cite{GS}, \cite{KR}, \cite{vanB}, \cite{NR}, \cite{BA});
  \item E-type/D-type tradition (\cite{EG}, \cite{NS}, \cite{HI}, \cite{EP});
  \item Modern type-theoretic approaches with dependent types (\cite{RA}, \cite{FT}, \cite{CR}, \cite{RT}, \cite{BD}, see also chapters one and three of this volume).
\end{itemize}
Our proposal belongs with the last group of modern type-theoretic approaches. The main technical contribution of our proposal consists in combining generalized quantifiers (\cite{MA}, \cite{LP}, \cite{BC}) with dependent types (\cite{MLP}, \cite{RA}, \cite{MM}). Empirically, our system allows a uniform treatment of all types of unbound anaphora, with the anaphoric (dynamic) effects falling out naturally as a consequence of having generalized quantification on dependent types.

The paper is organized as follows. Section 2 introduces informally the main features of our proposal. In this section we also describe our process of English-to-formal language translation. Section 3 shows how to interpret a range of anaphoric data in our system (maximal anaphora to quantifiers, quantificational subordination, cumulative and branching continuations, and `donkey anaphora'). Finally, sections 4 and 5 define the syntax and semantics of the system.

\section{Main features of the system}
The main elements of our system are
  \begin{itemize}
    \item Context and type dependency
    \item Many-typed (many-sorted) analysis
    \item Generalized quantifiers on dependent types
    \item Dynamic extensions of contexts
  \end{itemize}
  The discussion of the dynamic extensions of contexts is left for the next section.

 \subsection{Context and type dependency}

  The approaches adopted within the modern type-theoretic framework have been either proof-theoretic, where proof is a central semantic concept: \cite{RA}, \cite{L12a}, \cite{BD}, \cite{RT} (see also chapters one and two of this volume), or involved a combination of proof-theoretic and model-theoretic elements:  \cite{FT}, \cite{CR} (see also chapter three of this volume). By contrast to the existing proposals, our approach is model-theoretic with truth and reference being basic concepts (and no proofs). The two key type-theoretic features in our system are: context and type dependency.

 \subsubsection{Types, dependent types and their interpretation}

The variables of our system are always typed.
\begin{itemize}
  \item We write $x:X$ to denote that the variable $x$ is of type $X$ and refer to this as a type specification of the variable $x$.
  \item Types, in our system, are interpreted as sets. We write the interpretation of the type $X$ as $\|X\|$.
\end{itemize}
Types can depend on variables of other types.
\begin{itemize}
  \item If we already have a type specification $x:X$, then we can also have type $Y(x)$ depending on the variable $x$ and we can declare a variable $y$ of type $Y$ by stating $y:Y(x)$.
  \item The fact that $Y$ depends on $X$ is modeled as a function (projection)\\ $\pi: \|Y\| \rightarrow \|X\|.$
\end{itemize}
One example of such a dependence of types is that if $m$ is a variable of the type of months $M$, there is a type $D(m)$ of the days in that month

\vspace {4mm}

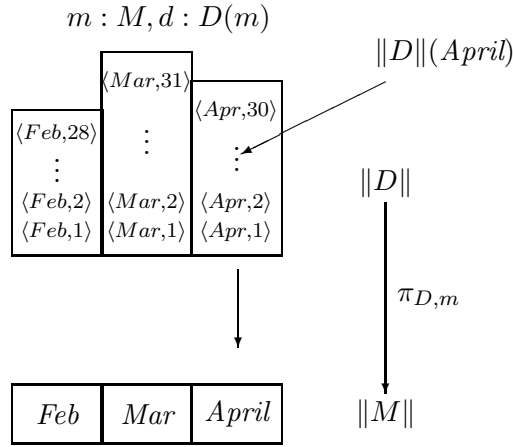
\begin{center}
 $m:M, d:D(m)$
\end{center}

\begin{center} \xext=1500 \yext=1300
\begin{picture}(\xext,\yext)(\xoff,\yoff)

  \put(210,0){\framebox(310,200){\em Feb}}
  \put(520,0){\framebox(310,200){\em Mar}}
  \put(830,0){\framebox(310,200){\em April}}

  \put(210,650){\framebox(310,500){}}
  \put(210,680){\makebox(310,100){\mbox{$_{\lk Feb,1\rk}$}}}
  \put(210,780){\makebox(310,100){\mbox{$_{\lk Feb,2\rk}$}}}
  \put(210,920){\makebox(310,100){\vdots}}
  \put(210,1020){\makebox(310,100){\mbox{$_{\lk Feb,28\rk}$}}}

  \put(520,650){\framebox(310,700){}}
  \put(520,680){\makebox(310,100){\mbox{$_{\lk Mar,1\rk}$}}}
  \put(520,780){\makebox(310,100){\mbox{$_{\lk Mar,2\rk}$}}}
  \put(520,1020){\makebox(310,100){\vdots}}
  \put(520,1200){\makebox(310,100){\mbox{$_{\lk Mar,31\rk}$}}}

  \put(830,650){\framebox(310,600){}}
  \put(830,680){\makebox(310,100){\mbox{$_{\lk Apr,1\rk}$}}}
  \put(830,780){\makebox(310,100){\mbox{$_{\lk Apr,2\rk}$}}}
  \put(830,960){\makebox(310,100){\vdots}}
  \put(830,1100){\makebox(310,100){\mbox{$_{\lk Apr,30\rk}$}}}

  \put(1550,1300){\makebox(310,100){\mbox{$\|D\|$({\em April})}}}
  \put(1500,1250){\vector(-2,-1){500}}

\put(990,600){\vector(0,-1){270}}

\putmorphism(1500,900)(0,-1)[\|D\|`\|M\|`\pi_{D,m}]{800}{1}r

 \end{picture}
\end{center}

\vspace{4mm}

\noindent If we interpret type $M$ as a set $\| M\|$ of months, then we can interpret type $D$ as a set of the days of the months in $\|M\|$, i.e. as a set of pairs
\[ \| D\| =\{ \lk a,k\rk  : k \; {\rm is\; (the\; number\; of) \; a\; day\; in\; month}\; a \}\]
equipped with the projection $\pi: \|D\| \rightarrow \|M\|$. The particular sets $\|D\|(a)$ of the days of the month $a$ can be recovered as the {\em fibers} of this projection (the preimages of $\{a\}$ under $\pi$)
\[\|D\|(a) = \{d \in \|D\|: \pi(d)=a\}.\]

\subsubsection{Contexts and their interpretation}

In type-theoretic settings, we can have a sequence of type specifications of the (individual) variables

\vspace{2mm}

\[\Gamma = x: X, y: Y(x), z: Z(x,y), t: T(x), u:U, \ldots \]
We adopt the convention that the variables the types depend on are always explicitly written down in specifications. Thus type $Y$ depends on the variable $x$; type $Z$, on the variables $x$ and $y$; type $T$, just on the variable $x$; and type $U$ is an example of a {\em constant type}, i.e. it does not depend on any variables. {\em Context} for us is a partially ordered sequence of type specifications of the (individual) variables such that the declaration of a variable $x$ (of type $X$) precedes the declaration of a variable $y$ (of type $Y$) if the type $Y$ depends on the variable $x$.

Contexts give rise to dependence graphs. A {\em dependence graph} for the context $\Gamma$ is a graph that has types occurring in $\Gamma$ as vertices, and, for every variable specification $x:X(\ldots)$ and type $Y(\ldots, x,\ldots)$ that depends on $x$ in $\Gamma$, it has an edge

\begin{center} \xext=200 \yext=700
\begin{picture}(\xext,\yext)(\xoff,\yoff)
\putmorphism(0,600)(0,-1)[Y`X`\pi_{Y,x}]{500}{1}r
 \end{picture}
\end{center}
The corresponding semantic notion is that of a dependence diagram. The {\em dependence diagram} for the context $\Gamma$ associates to every type $X$ in $\Gamma$ a set $\|X\|$, and to every edge $\pi_{Y,x} : Y \ra X$, a function $\|\pi_{Y,x}\| : \|Y\| \ra \|X\|$, so that whenever we have a triangle of edges (as on the left), the corresponding triangle of functions commutes (i.e. $\|\pi_{Z,x}\| = \|\pi_{Y,x}\| \circ \|\pi_{Z,y}\|$)

\vspace{2mm}

\begin{center} \xext=2000 \yext=650
\begin{picture}(\xext,\yext)(\xoff,\yoff)

  \settriparms[1`1`1;300]
  \putDtriangle(0,0)[Z`Y`X;\pi_{Z,x}`\pi_{Z,y}`\pi_{Y,x}]

  \settriparms[1`1`1;300]
  \putDtriangle(1500,0)[\|Z\|`\|Y\|`\|X\|;\|\pi_{Z,x}\|`\|\pi_{Z,y}\|`\|\pi_{Y,x}\|]
 \end{picture}
\end{center}

\vspace{2mm}

\noindent We say that, for $a \in \|X\|, b \in \|Y\|, c \in \|Z\|$, a triple $\langle a, b, c \rangle$ is {\em compatible} iff

\[ \pi_{Y,x}(b) = a, \hskip 5 mm \pi_{Z,y}(c) = b, \hskip 5 mm  \pi_{Z,x}(c) = a. \]

The interpretation of the context $\Gamma$, the {\em parameter space} $\| \Gamma\|$, is a set of compatible $n$-tuples of the elements of the sets corresponding to the types involved (compatible wrt all projections)

\[\| \Gamma\| = \| x_1: X_1, \ldots, x_n : X_n(\lk x_i\rk_{i\in J_n}) \| =  \{  \lk \bar{x}_1, \ldots, \bar{x}_n \rk : \bar{x}_i\in \|X_i\|,\; {\rm and}\; \] \[ \hskip 9 cm \|\pi_{X_{i'},x_i}\|(\bar{x}_{i'})= \bar{x}_i  \} \]

\subsection{Many-typed (many-sorted) analysis}

Like in the classical Montague-style approach, we have generalized quantifiers in our system. But in the spirit of the modern type-theoretic framework we adopt a many-typed analysis (in place of a standard single-sorted analysis). Such richer type structures have been also extensively applied to studies of lexical phenomena such as selection restriction or coersions (\cite{asher}, \cite{L12a}, \cite{L12b}, \cite{retore}, see also chapters two, five and six of this volume).

\subsubsection{Montague-style semantics}
Standard Montague-style semantics is single-sorted in the sense that it includes one type \textbf{e} of all entities (strictly speaking, it has two basic types: type \textbf{e} and type  \textbf{t} of truth values, and a recursive definition of functional types); quantifiers and predicates are interpreted over the universe of all entities $E$.

On the Montague-style analysis, quantifier phrases, e.g. \textit{every man} or \textit{some woman}, are interpreted as sets of subsets of $E$
 \[\|every \; man\| = \{X \subseteq E: \|man\| \subseteq X\}.\]
 \[\|some \; woman\| = \{X \subseteq E: \|woman\| \cap X \neq \emptyset\}.\]
 On this standard analysis, an element of the denotation of a quantifier phrase like \textit{every man} or \textit{some woman} (i.e. a subset of the universe, $X \subseteq E$) will contain besides men or women all sorts of entities (children, books, etc). To have elements from which such extra entities are removed, Barwise and Cooper define notions such as `witness set' (see \cite{BC}, \cite{SA}). Quantifier phrases are interpreted this way to ensure that predicates are unambiguous. On the Montague-style analysis, a predicate like \textit{love} denotes a single \textit{love}-relation, whether relating men to women, children to mothers, etc.

\subsubsection{Polymorphic interpretation of quantifiers and predicates}

Our analysis is many-sorted in the sense that it includes many basic types, and so we have a polymorphic interpretation of quantifiers and predicates.

A generalized quantifier associates to every set $Z$ a subset of the power set of $Z$
 \[\|Q\|(Z)\subseteq \cP(Z)\]
Quantifier phrases, e.g. \textit{every man} or \textit{some woman}, are interpreted as
 \[ \|\forall_{m:Man}\| = \{\|Man\|\}\]
 \[ \|\exists_{w:Woman}\| = \{X\subseteq \|Woman\| : \hspace{1mm} X \neq \emptyset \}\]
  \textit{Every man} denotes a singleton set whose only element is the entire set of men (given in the context); \textit{some woman} denotes the set of all non-empty subsets of the set of women. As an element of the denotation of a quantifier phrase \textit{every man} or \textit{some woman} is homogeneous (containing men or women only), we do not need to consider notions such as `witness set'. As a consequence of our many-typed analysis, predicates are also defined polymorphically. If we have a predicate $P$ defined in a context $\Gamma$
\[  x_1:X_1, \ldots, x_n:X_n(\lk x_i\rk_{i\in J_n}) \vdash  P(\vec{x})\]
then, for any interpretation of the context $\|\Gamma\|$, the predicate is interpreted as a subset of its parameter space, i.e.  $\|P\|\subseteq \|\Gamma\|$.

\subsection{Generalized quantifiers on dependent types}

The interpretation of quantifier phrases is further extended into the interpretation of (generalized) quantifier prefixes.

\subsubsection{Combining quantifier phrases - chains of quantifiers}

Multi-quantifier sentences such as \textit{Every man loves a woman} or \textit{Last year two scientists wrote five papers} have been known to be ambiguous with different readings corresponding to how various quantifiers are semantically related in the sentence. To account for the readings available for such multi-quantifier sentences, we raise quantifier phrases to the front of a sentence to form (generalized) quantifier prefixes - {\em chains of quantifiers}. Chains of quantifiers are built from quantifier phrases using three chain-constructors: pack-formation rule $(?,\ldots,?)$, sequential composition $?|?$, and parallel composition $\frac{\hskip 2mm ?\hskip 2mm}{?}$. More precisely, quantifier phrases can be grouped together to form packs of quantifiers (one-element packs are considered quantifier phrases); (pre-)chains are then built from packs via the chain-constructors of sequential and parallel composition. The semantical operations that correspond to the three chain-constructors allow us to capture in a compositional manner cumulative, scope-dependent and branching readings

\vspace{4mm}

\[
\begin{array}{|c|c|c|c|} \hline

  \textit{chain constructors} & \textit{semantical operations} \\ \hline

  \textit{pack formation rule $(?,\ldots,?)$} & \textit{cumulation} \\ \hline

  \textit{sequential composition $?|?$} & \textit{iteration} \\ \hline

  \textit{parallel composition $\frac{\hskip 2mm ?\hskip 2mm}{?}$} & \textit{branching} \\ \hline
\end{array}
\]

\vspace{4mm}

To illustrate the working of the chain constructors and their corresponding semantical operations, we will first use a familiar example. \textit{Every man loves a woman} can be understood to mean that each of the men loves a potentially different woman. To capture this reading

\begin{itemize}
  \item a {\em sequential composition constructor} $?|?$ is used to produce a multi-quantifier prefix (chain of quantifiers): $\forall_{m:M}|\exists_{w:W}$;
  \item the corresponding semantical operation of {\em iteration} is defined as follows
\small \[\|\forall_{m:M}|\exists_{w:W}\| = \{R \subseteq \|M \|\times \|W\|: \hskip 6cm \]
\[ \{a \in \|M\|:\{b \in \|W\|: \langle a, b\rangle \in R\}\in \|\exists_{w:W}\| \}\in \|\forall_{m:M}\|\}.\hskip 5cm\]
\end{itemize}
The chain $\forall_{m:M}|\exists_{w:W}$ denotes a set of relations such that the set of men such that each man is in this relation to at least one woman is the set of all men. Obviously, the iteration rule gives the same result as the standard nesting of quantifiers in first-order logic. The idea of chain-constructors and the corresponding semantical operations builds on Mostowski's notion of quantifier (\cite{MA}) further generalized by Lindstr\"{o}m to a so-called polyadic quantifier (\cite{LP}). (See \cite{BZ}, compare also \cite{KE87}, \cite{BvJ}, \cite{KE92}, \cite{KE93}, \cite{WD}). A quantifier phrase like $\exists_{w:Woman}$ can be thought of as a one-place (monadic) quantifier and has as denotation a set of sets. A chain of quantifiers like $\forall_{m:M}|\exists_{w:W}$ can be thought of as a single two-place (polyadic) quantifier and has as denotation a set of binary relations.

Consider now a cumulative example. \textit{Last year three scientists wrote five papers} allows a reading saying that each of the three scientists wrote at least one of the five papers, and each of the five papers was written by at least one of the two scientists. To capture the cumulative reading
\begin{itemize}
  \item a {\em pack formation rule} $(?,\ldots,?)$ is used to produce a multi-quantifier prefix (pack of quantifiers): $(Three_{s:S}, Five_{p:P})$.
  \item the corresponding semantical operation of {\em cumulation} is defined as follows
   \[\|(Three_{s:S}, Five_{p:P})\| = \hskip 9 cm \]
   \[ =\{R \subseteq \|S \|\times \|A\|: \pi_1 (R)\in \|Three_{s:S}\| \; {\rm and}\; \pi_2(R)\in \|Five_{p:P}\|\} \hskip 8cm \]
   where $\pi_i$ is the $i$-th projection from the product.
\end{itemize}
Yet another reading is a branching reading where each of the three scientists wrote the same set of five papers. To capture this reading
\begin{itemize}
  \item  a {\em parallel composition constructor} $\frac{\hskip 2mm ?\hskip 2mm}{?}$ is used to produce a multi-quantifier prefix (chain of quantifiers): $\frac{Three_{s:S}}{Five_{p:P}}$.
  \item the corresponding semantical operation of {\em branching} is defined as follows
        \[ \|\frac{Three_{s:S}}{Five_{p:P}}\| = \{ A\times B \; :\; A \in \|Three_{s:S}\| \; {\rm and}\; B\in \|Five_{p:P}\|\} \hskip 8cm \]
\end{itemize}

\subsubsection{Combining generalized quantifiers with dependent types}

The three chain-constructors and the corresponding semantical operations are further extended to dependent types. To use an example of the iteration operation, we have
\small \[\|\forall_{m : M}|\exists_{w_D : W_D(m)}\| = \{R \subseteq \|W_D\|:\{a \in \|M\|: \hskip 6cm \]
\[ \{b \in \|W_D\|(a): \langle a, b\rangle \in R\}\in \|\exists_{w_D: W_D(m)}\|(\|W_D\|(a)) \}\in \|\forall_{m:M}\|\}. \hskip 6cm \]

\vspace{2mm}

\normalsize
\noindent The chain $\forall_{m:M}|\exists_{w_D:W_D(m)}$ denotes a set of relations such that the set of men such that each man is in this relation to \textbf{at least one woman in the corresponding fiber of women} is the set of all men. By extending chains of quantifiers to dependent types, our system introduces {\em quantification over fibers} - in the example used, existential quantification over fibers of women $\|W_D\|(a)$

\vspace{4mm}

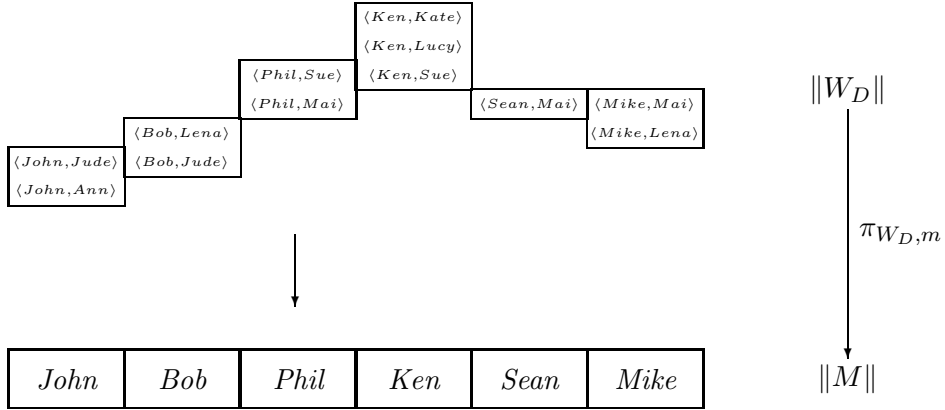
\begin{center} \xext=3000 \yext=1400
\begin{picture}(\xext,\yext)(\xoff,\yoff)

  \put(0,0){\framebox(400,200){\em John}}
  \put(400,0){\framebox(400,200){\em Bob}}
  \put(800,0){\framebox(400,200){\em Phil}}
  \put(1200,0){\framebox(400,200){\em Ken}}
  \put(1600,0){\framebox(400,200){\em Sean}}
  \put(2000,0){\framebox(400,200){\em Mike}}

  \put(0,700){\framebox(400,200){}}
  \put(0,700){\makebox(400,100){\mbox{\tiny $_{\lk John,Ann\rk}$}}}
  \put(0,800){\makebox(400,100){\mbox{\tiny $_{\lk John,Jude\rk}$}}}

  \put(400,800){\framebox(400,200){}}
  \put(400,800){\makebox(400,100){\mbox{\tiny $_{\lk Bob,Jude\rk}$}}}
  \put(400,900){\makebox(400,100){\mbox{\tiny $_{\lk Bob,Lena\rk}$}}}

  \put(800,1000){\framebox(400,200){}}
  \put(800,1000){\makebox(400,100){\mbox{\tiny $_{\lk Phil,Mai\rk}$}}}
  \put(800,1100){\makebox(400,100){\mbox{\tiny $_{\lk Phil,Sue\rk}$}}}

  \put(1200,1100){\framebox(400,300){}}
  \put(1200,1100){\makebox(400,100){\mbox{\tiny $_{\lk Ken,Sue\rk}$}}}
  \put(1200,1200){\makebox(400,100){\mbox{\tiny $_{\lk Ken,Lucy\rk}$}}}
  \put(1200,1300){\makebox(400,100){\mbox{\tiny $_{\lk Ken,Kate\rk}$}}}

  \put(1600,1000){\framebox(400,100){}}
  \put(1600,1000){\makebox(400,100){\mbox{\tiny $_{\lk Sean,Mai\rk}$}}}

  \put(2000,900){\framebox(400,200){}}
  \put(2000,900){\makebox(400,100){\mbox{\tiny $_{\lk Mike,Lena\rk}$}}}
  \put(2000,1000){\makebox(400,100){\mbox{\tiny $_{\lk Mike,Mai\rk}$}}}

\put(990,600){\vector(0,-1){250}}

\putmorphism(2900,1100)(0,-1)[\|W_D\|`\|M\|`\pi_{W_D,m}]{1000}{1}r

 \end{picture}
\end{center}

\vskip 2mm

\noindent In this sense, fibers are considered 1st class citizens of our semantics, i.e. our system allows for quantification over fibers on a par with quantification over any other type.

Note that in a system with generalized quantification extended to dependent types, chains of quantifiers are composed out of pre-chains

\begin{itemize}
  \item $Ch_{\vec{y}:\vec{Y}(\vec{x})}$  denotes a {\em pre-chain} with {\em binding variables} $\vec{y}$ and {\em indexing variables} $\vec{x}$.
  \item Chains of quantifiers are pre-chains in which all indexing variables are bound.
\end{itemize}
In order to make sure that a pre-chain can be turned into a chain, we impose a global restriction on variables that each occurrence of an indexing variable in $Ch$ be preceded by a binding occurrence of that variable in $Ch$. Below we give examples of both correct and incorrect pre-chains

\vskip 2mm

Correct pre-chain
\[\frac{{Q_2}_{y:Y(x)}|{Q_3}_{z:Z(x,y)}}{{Q_4}_{u:U}} \]
The pre-chain above can be turned into a chain, e.g. by prefixing ${Q_1}_{x:X}$ and binding indexing occurrences of $x$
\[\frac{{Q_1}_{x:X}|{Q_2}_{y:Y(x)}|{Q_3}_{z:Z(x,y)}}{{Q_4}_{u:U}} \]

\vskip 2mm

Incorrect pre-chain
\[\frac{{Q_1}_{y:Y(x)}|{Q_2}_{x:X(y,z)}}{{Q_3}_{u:U}} \]
The pre-chain above is incorrect, as the occurrence of an indexing variable $x$ is followed by the binding occurrence of that variable and so cannot
get bound.

\subsection{English-to-formal language translation}\label{1}

Our English-to-formal language translation process consists of two steps (i) {\em representation} and (ii) {\em disambiguation}.

\noindent {\em Representation.} The syntax of the representation language - for the English fragment considered in this paper - is as follows

\vskip 2mm

\noindent $S \rightarrow Prd^{n}(QP_{1}, \ldots , QP_{n})$;\\
$MCN \rightarrow Prd^{n}(QP_{1}, \ldots, \; CN \; ,\ldots , QP_{n});$\\
$MCN \rightarrow CN;$\\
$QP \rightarrow Det \; MCN$;\\
$Det \rightarrow every, most, three, \ldots$;\\
$CN \rightarrow man, woman,\ldots$;\\
$Prd^{n} \rightarrow enter, love, \ldots$

\vskip 2mm

\noindent In the Montague-style semantics, common nouns ($CN$) are interpreted as predicates (expressions of type $e \rightarrow t$). In our type-theoretic setting, $CNs$ are interpreted as types; modified common nouns ($MCNs$, to be discussed below), as $\ast$-sentences determining some (possibly dependent) types, and predicates are interpreted over the types.

\noindent {\em Disambiguation.} Sentences of English, contrary to sentences of our formal language, are often ambiguous. Hence one sentence representation
can be associated with more than one sentence in our formal language. The second step thus involves disambiguation. We take quantifier phrases of a given representation, e.g.
\begin{center}
$P(Q_1X_1,Q_2X_2,Q_3X_3)$
\end{center}
and organize them into all possible chains of quantifiers in suitable contexts with some restrictions imposed on particular quantifiers concerning the places in prefixes at which they can occur (a detailed elaboration of the disambiguation process is left for another place)
\begin{center}
$\frac{Q_1 x_1: X_1|Q_2 x_2: X_2}{Q_3 x_3: X_3} \;P(x_1, x_2, x_3).$
\end{center}

\section{Dynamic extensions of contexts}

Our interpretational architecture is two-dimensional. The two dimensions to the meaning of a sentence in our system are: the truth value of a sentence and the dynamic effects introduced by the sentence.

A sentence with a chain of quantifiers $Ch_{\vec{y}:\vec{Y}}$ and predicate $P(\vec{y})$ is true iff the interpretation of the predicate (i.e. some set of compatible $n$-tuples) belongs to the interpretation of the chain (i.e. some family of sets of compatible $n$-tuples), i.e. iff
\[\|P\|(\|\vec{y}:\vec{Y}\|)\in \|Ch_{\vec{y}:\vec{Y}}\|.\]

A sentence with a chain of quantifiers also extends the context, i.e. it creates a new context out of the old one (in which it takes place) by adding some possibly dependent types; the anaphoric continuation then is interpreted in the newly obtained context
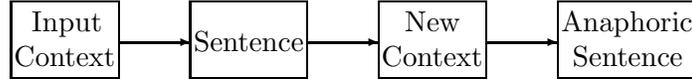
\begin{center}
\xext=2200 \yext=300
\begin{picture}(\xext,\yext)(\xoff,\yoff)

 \put(70,160){Input}
 \put(10,40){Context}
 \put(0,0){\framebox(370,260){ }}
  \put(370,120){\vector(1,0){250}}

 \put(630,90){Sentence}
 \put(620,0){\framebox(400,260){ }}
  \put(1020,120){\vector(1,0){250}}

  \put(1360,160){New}
 \put(1280,40){Context}
 \put(1270,0){\framebox(370,260){ }}
  \put(1640,120){\vector(1,0){250}}

  \put(1900,160){Anaphoric}
 \put(1940,40){Sentence}
 \put(1890,0){\framebox(470,260){ }}

\end{picture}
\end{center}
For the purpose of modeling the dynamic extensions of context, we introduce a new type constructor $\T$. For the interpretation of the types from the extended context, we define a new algorithm. We now show how to interpret a range of anaphoric data in our system: maximal anaphora to quantifiers, quantificational subordination (including iterated examples), cumulative and branching continuations, and `donkey anaphora' (including iterated `donkey sentences').

 \subsection{Maximal anaphora to quantifiers}

Let us first consider an example in (1)
\begin{enumerate}[(1)]
\item\label{example1} Most kids entered. They looked happy.
\end{enumerate}
As already mentioned in Section \ref{4}, the observation is that the anaphoric pronoun \textit{they} in the second sentence refers to the entire set of kids who entered. Thus the first sentence must introduce the set of all kids who entered.

We start with \textsc{Input Context}
\[\Gamma:=\; k: Kid \]
Sentence $\varphi:=$ \textit{Most kids entered} translates into a sentence with a chain of quantifiers in the \textsc{Input Context} $\Gamma$
\[\Gamma \vdash Most_{k : K} Enter(k),\]
and creates \textsc{New Context} by adding a new variable specification on a newly formed type
\[\Gamma_\varphi:=\; k: Kid, \;\ t_{\varphi, Most_k}: \T_{\varphi, Most_{k : K}}\]
$\textit{Anaphoric continuation}: =$ \textit{They looked happy} is now interpreted in the \textsc{New Context} $\Gamma_\varphi$
\[\Gamma_\varphi \vdash \forall_{t_{\varphi, Most_k}: \T_{\varphi, Most_{k : K}}} \,\,  Happy(t_{\varphi, Most_k}).\]
 We follow here E-type/D-type tradition (\cite{EG}, \cite{NS}, \cite{HI}, \cite{EP}) in assuming that unbound anaphoric pronouns are subject to a maximality constraint, i.e. by default they are treated as universal quantifiers; context is used as a medium supplying possibly dependent types as their potential quantificational domains.

The interpretation of the new type from the extended context is defined by our procedure as
\[ \|\T_{\varphi, Most_{k : K}}\|:=\|Enter\|\]

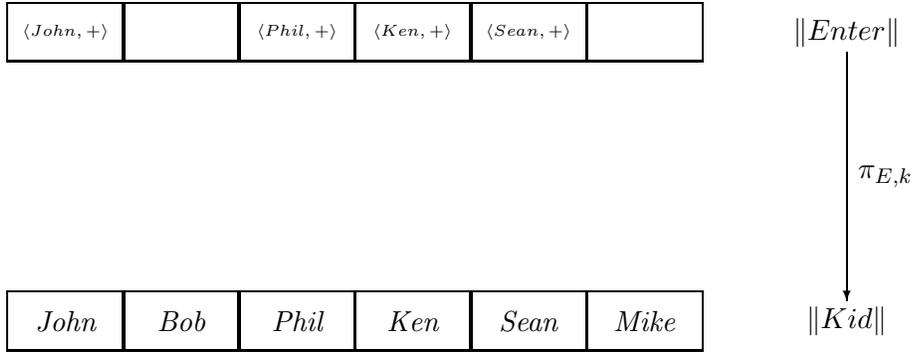
\begin{center} \xext=3000 \yext=1400
\begin{picture}(\xext,\yext)(\xoff,\yoff)

  \put(0,0){\framebox(400,200){\em John}}
  \put(400,0){\framebox(400,200){\em Bob}}
  \put(800,0){\framebox(400,200){\em Phil}}
  \put(1200,0){\framebox(400,200){\em Ken}}
  \put(1600,0){\framebox(400,200){\em Sean}}
  \put(2000,0){\framebox(400,200){\em Mike}}

  \put(0,1000){\framebox(400,200){\tiny $\lk John, +\rk$}}
  \put(400,1000){\framebox(400,200){}}
  \put(800,1000){\framebox(400,200){\tiny $\lk Phil, +\rk$}}
  \put(1200,1000){\framebox(400,200){\tiny $\lk Ken, +\rk$}}
  \put(1600,1000){\framebox(400,200){\tiny $\lk Sean, +\rk$}}
  \put(2000,1000){\framebox(400,200){}}

\putmorphism(2900,1100)(0,-1)[\|Enter\|`\|Kid\|`\pi_{E,k}]{1000}{1}r

 \end{picture}
\end{center}
  Thus on our analysis the pronoun \textit{they} in the second sentence quantifies universally over the set $\|Enter\|$, yielding the correct truth-conditions for the anaphoric continuation \textit{Every kid who entered looked happy}.

\subsection{Quantificational subordination}

Consider now a case of quantificational subordination (to better illustrate the full benefits of our interpretational algorithm, we will use a more difficult variant of the example introduced in Section \ref{4})
\begin{enumerate}[(2)]
\item Most men love two women. They (each) kiss them.
\end{enumerate}
The first sentence in (2) (on the interpretation where \textit{two women} depends on \textit{most men}) is understood to mean that most men are such that they each love a potentially different set of two women. The way to understand the second sentence in (2) is that every man who loves two women kisses the women he loves rather than those loved by someone else. Thus, intuitively, the first sentence in (2) must deliver  a dependency between each of the men and the women they love.

We start with \textsc{Input Context}
\[\Gamma:= \; m: Man, w: Woman\]
Sentence $\varphi :=$ \textit{Most men love two women} translates into a sentence with a chain of quantifiers in the \textsc{Input Context}  $\Gamma$
\[\Gamma \vdash Most_{m:M}| Two_{w:W} Love(m, w).\]
and creates \textsc{New Context} by adding new variable specifications on two newly formed types
\small \[\Gamma_\varphi := \; m: Man, w: Woman, t_{\varphi, Most_m}: \T_{\varphi, Most_{m:M}}; \  t_{\varphi, Two_w} : \T_{\varphi, Two_{w:W}}(t_{\varphi, Most_m})\]
\normalsize $\textit{Anaphoric continuation}: =$ \textit{They (each) kiss them} is interpreted in the \textsc{New Context} $\Gamma_\varphi$
\[ \Gamma_\varphi \vdash \forall_{t_{\varphi, Most_m}: \T_{\varphi, Most_{m:M}}} | \forall_{t_{\varphi, Two_w} : \T_{\varphi, Two_{w:W}}(t_{\varphi, Most_m})} Kiss (t_{\varphi, Most_m}, t_{\varphi, Two_w}). \hskip 8cm\]

The interpretations of the types from the extended context are defined in a two-step procedure.

\noindent \textbf{Step 1.} We define fibers of new types (by inverse induction from chains down to quantifier phrases).

\noindent Basic step. For the whole chain $Ch = Most_{m:M}|Two_{w:W}$ we put

\[\|\T_{\varphi, Most_{m:M}| Two_{w:W}}\|:=\|Love\| \]
i.e. we take the interpretation of $\T_{Ch}$ to be the denotation of the whole predicate $\|Love\|$.

\vskip 2mm

\noindent Inductive step.

\vskip 2mm

\noindent For $a\in \|M\|$,
\[\|\T_{\varphi, Two_{w:W}}\|(a) = \{b\in \|W\| \, : \, \lk a,b \rk\in \|Love\| \} \hskip 8cm\]

\[\|\T_{\varphi, Most_{m:M}}\| =  \{a \in \|M\| \, : \,  \{b\in \|W\| \, : \, \lk a,b \rk\in \|Love\| \} \in  \|Two_{w:W}\|\} \hskip 8cm\]

\vspace{2mm}

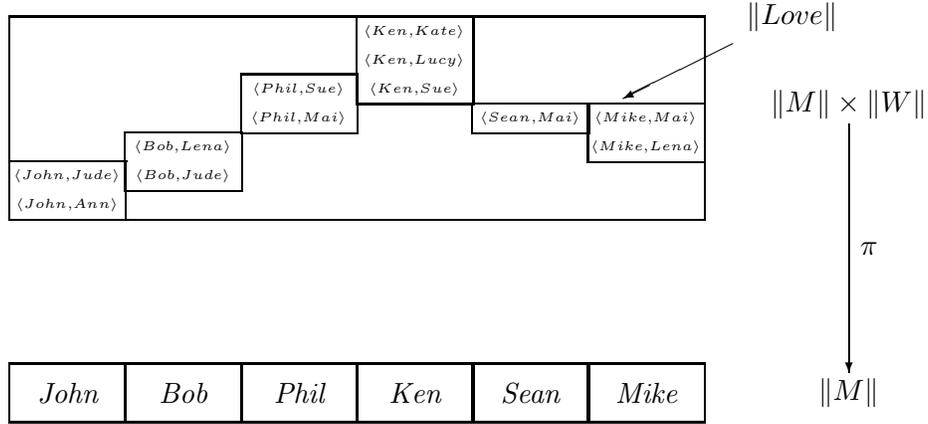
\begin{center} \xext=3200 \yext=1400
\begin{picture}(\xext,\yext)(\xoff,\yoff)
\put(0,700){\framebox(2400,700){}}

  \put(0,0){\framebox(400,200){\em John}}
  \put(400,0){\framebox(400,200){\em Bob}}
  \put(800,0){\framebox(400,200){\em Phil}}
  \put(1200,0){\framebox(400,200){\em Ken}}
  \put(1600,0){\framebox(400,200){\em Sean}}
  \put(2000,0){\framebox(400,200){\em Mike}}

  \put(0,700){\framebox(400,200){}}
  \put(0,700){\makebox(400,100){\mbox{\tiny $_{\lk John,Ann\rk}$}}}
  \put(0,800){\makebox(400,100){\mbox{\tiny $_{\lk John,Jude\rk}$}}}

  \put(400,800){\framebox(400,200){}}
  \put(400,800){\makebox(400,100){\mbox{\tiny $_{\lk Bob,Jude\rk}$}}}
  \put(400,900){\makebox(400,100){\mbox{\tiny $_{\lk Bob,Lena\rk}$}}}

  \put(800,1000){\framebox(400,200){}}
  \put(800,1000){\makebox(400,100){\mbox{\tiny $_{\lk Phil,Mai\rk}$}}}
  \put(800,1100){\makebox(400,100){\mbox{\tiny $_{\lk Phil,Sue\rk}$}}}

  \put(1200,1100){\framebox(400,300){}}
  \put(1200,1100){\makebox(400,100){\mbox{\tiny $_{\lk Ken,Sue\rk}$}}}
  \put(1200,1200){\makebox(400,100){\mbox{\tiny $_{\lk Ken,Lucy\rk}$}}}
  \put(1200,1300){\makebox(400,100){\mbox{\tiny $_{\lk Ken,Kate\rk}$}}}

  \put(1600,1000){\framebox(400,100){}}
  \put(1600,1000){\makebox(400,100){\mbox{\tiny $_{\lk Sean,Mai\rk}$}}}

  \put(2000,900){\framebox(400,200){}}
  \put(2000,900){\makebox(400,100){\mbox{\tiny $_{\lk Mike,Lena\rk}$}}}
  \put(2000,1000){\makebox(400,100){\mbox{\tiny $_{\lk Mike,Mai\rk}$}}}

\putmorphism(2900,1100)(0,-1)[\|M\|  \times \|W\|`\|M\|`\pi]{1000}{1}r

\put(2550,1350){\makebox(310,100){\mbox{$\|Love\|$}}}
\put(2500,1310){\vector(-2,-1){380}}

 \end{picture}
\end{center}

\vspace{2mm}

\noindent \textbf{Step 2.} We build dependent types from fibers.

\[ \|\T_{\varphi, Two_{w:W}}\| = \bigcup\{ \{ a\}\times \|\T_{\varphi, Two_{w:W}}\|(a) : a \in \|\T_{\varphi, Most_{m:M}}\|\} \hskip 8cm\]

\[ \|\T_{\varphi, Most_{m:M}}\| =  \{a \in \|M\| \, : \,  \{b\in \|W\| \, : \, \lk a,b \rk\in \|Love\| \} \in  \|Two_{w:W}\| \} \hskip 8cm\]

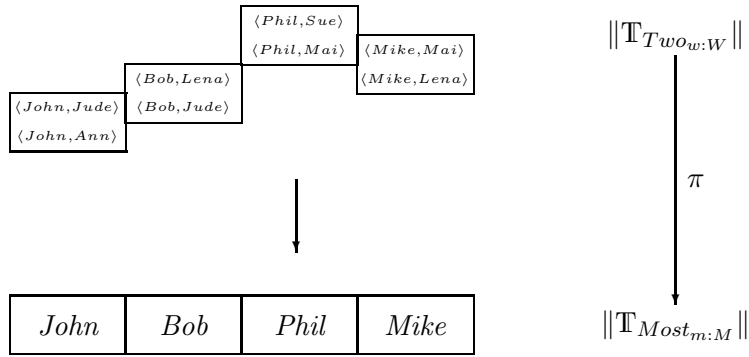
\begin{center} \xext=2500 \yext=1400
\begin{picture}(\xext,\yext)(\xoff,\yoff)

  \put(0,0){\framebox(400,200){\em John}}
  \put(400,0){\framebox(400,200){\em Bob}}
  \put(800,0){\framebox(400,200){\em Phil}}
  \put(1200,0){\framebox(400,200){\em Mike}}

  \put(0,700){\framebox(400,200){}}
  \put(0,700){\makebox(400,100){\mbox{\tiny $_{\lk John,Ann\rk}$}}}
  \put(0,800){\makebox(400,100){\mbox{\tiny $_{\lk John,Jude\rk}$}}}

  \put(400,800){\framebox(400,200){}}
  \put(400,800){\makebox(400,100){\mbox{\tiny $_{\lk Bob,Jude\rk}$}}}
  \put(400,900){\makebox(400,100){\mbox{\tiny $_{\lk Bob,Lena\rk}$}}}

  \put(800,1000){\framebox(400,200){}}
  \put(800,1000){\makebox(400,100){\mbox{\tiny $_{\lk Phil,Mai\rk}$}}}
  \put(800,1100){\makebox(400,100){\mbox{\tiny $_{\lk Phil,Sue\rk}$}}}

  \put(1200,900){\framebox(400,200){}}
  \put(1200,900){\makebox(400,100){\mbox{\tiny $_{\lk Mike,Lena\rk}$}}}
  \put(1200,1000){\makebox(400,100){\mbox{\tiny $_{\lk Mike,Mai\rk}$}}}

\put(990,600){\vector(0,-1){250}}

\putmorphism(2300,1100)(0,-1)[\|\T_{Two_{w:W}}\|`\|\T_{Most_{m:M}}\|`\pi]{1000}{1}r

 \end{picture}
\end{center}

Thus the context gets extended by
\begin{itemize}
\item the type interpreted as $\|\T_{Most_{m:M}}\|$, i.e.\ the set of men who love two women;
\item the dependent type interpreted for $a \in \|\T_{Most_{m:M}}\|$ as $\|\T_{Two_{w:W}}\|(a)$, i.e.\ the set of women loved by the man $a$.
\end{itemize}
The two unbound anaphoric pronouns $they_m$ and $them_w$ in the second sentence of (2) quantify universally over the respective interpretations, yielding the correct truth conditions \textit{Every man who loves two women kisses every woman he loves}. Note that the anaphoric continuation in this example crucially involves (universal) quantification over fibers of the women loved, $\|\T_{\varphi, Two_{w:W}}\|(a)$.

\subsection{Cumulative and branching continuations}

Our system defines dynamic extensions of contexts and their interpretation also for cumulative and branching continuations. Consider examples in (3a) and (3b)
\begin{enumerate}[(3a)]
\item\label{example3a} Last year three scientists wrote (a total of) five articles (between them). They presented them at major conferences.
\end{enumerate}
\begin{enumerate}[(3b)]
\item\label{example3b} Last year three scientists (each) wrote (the same) five articles. They presented them at major conferences.
\end{enumerate}
As already discussed in Section \ref{4}, the dynamics of the first sentence in (3a) and (3b) can deliver some (respectively: cumulative or branching) internal relation between the types corresponding to \textit{three scientists} and \textit{five articles} that can be elaborated upon in the anaphoric continuation.

Consider first the cumulative example. The anaphoric continuation in (3a) can be interpreted in what Krifka calls a `correspondence' fashion (see \cite{KM}). For example, John wrote one article, co-authored one more with Bob, who co-authored one more with Ken who wrote two more articles by himself, and the scientists that cooperated in writing one or more articles also cooperated in presenting these (and no other) articles at major conferences.

We start with \textsc{Input Context}
\[ \Gamma :=\; s: Scientist, a: Article\]
Sentence $\varphi:=$ \textit{Three scientists wrote a total of five articles (between them)} translates into a sentence with a chain of quantifiers in the \textsc{Input Context} $\Gamma$
\[\Gamma \vdash (Three_{s: Scientist}, Five_{a: Article})\,\, Write(s, a),\]
and creates \textsc{New Context} by adding a new variable specification on a newly formed type
\[\Gamma_\varphi:=\; s: Scientist, a: Article, t_{\varphi, (Three_s, Five_a)}: \T_{\varphi, (Three_{s: S}; \ Five_{a: A})}\]
$\textit{Anaphoric continuation}: =$ \textit{They presented them at major conferences} is interpreted in the \textsc{New Context} $\Gamma_\varphi$
\[ \Gamma_\varphi \vdash \forall_{t_{\varphi, (Three_s, Five_a)}} Present (t_{\varphi, (Three_s, Five_a)}).\]

The interpretation of the new type from the extended context is defined by our procedure as
\[\|\T_{\varphi, (Three_{s: S}, Five_{a: A})}\| = \|Write\|. \]
\begin{center} \xext=1800 \yext=950
\begin{picture}(\xext,\yext)(\xoff,\yoff)
\put(0,100){\framebox(1800,500){}}

  \put(0,100){\framebox(600,200){}}
  \put(0,100){\makebox(600,100){\mbox{$_{\lk John,article 1\rk}$}}}
  \put(0,200){\makebox(600,100){\mbox{$_{\lk John,article 2\rk}$}}}

  \put(600,200){\framebox(600,200){}}
  \put(600,200){\makebox(600,100){\mbox{$_{\lk Bob, article 2\rk}$}}}
  \put(600,300){\makebox(600,100){\mbox{$_{\lk Bob, article 3\rk}$}}}

  \put(1200,300){\framebox(600,300){}}
  \put(1200,300){\makebox(600,100){\mbox{$_{\lk Ken, article 3\rk}$}}}
  \put(1200,400){\makebox(600,100){\mbox{$_{\lk Ken, article 4\rk}$}}}
  \put(1200,500){\makebox(600,100){\mbox{$_{\lk Ken, article 5\rk}$}}}

\put(2000,750){\makebox(310,100){\mbox{$\|Write\|\subseteq \| S\| \times \|A\|$}}}
\put(2000,710){\vector(-2,-1){380}}

 \end{picture}
\end{center}
The anaphoric continuation quantifies universally over the respective interpretation (i.e.\ a set of $\lk scientist, article \rk$ pairs such that the scientist wrote the article), yielding the desired truth-conditions \textit{The respective scientists cooperated in presenting at major conferences the respective articles that they cooperated in writing}.

Consider now the branching example. The way to understand the anaphoric continuation is that the three scientists - say, John, Bob and Ken - co-authored all of the five articles, and all of the scientists involved presented at major conferences all of the articles involved.

We start with \textsc{Input Context}
\[\Gamma:=\; s: Scientist, a: Article\]
Sentence $\varphi:=$ \textit{Three scientists (each) wrote (the same) five articles} translates into a sentence with a chain of quantifiers in the \textsc{Input Context} $\Gamma$
\[\Gamma \vdash\frac{Three_{s: S}}{Five_{a: A}} \,\, Write(s, a),\]
and creates \textsc{New Context} by adding two new variable specification on two newly formed constant types
\[\Gamma_\varphi:=\; s: Scientist, a: Article, t_{\varphi, Three_s}: \T_{\varphi, Three_{s: S}}; \ t_{\varphi, Five_a}: \T_{\varphi, Five_{a: A}}\]
$\textit{Anaphoric continuation}: =$ \textit{They presented them at major conferences} is interpreted in the \textsc{New Context} $\Gamma_\varphi$
\[\Gamma_\varphi \vdash \frac{\forall_{t_{\varphi, Three_s}}} {\forall_{t_{\varphi, Five_a}}} Present (t_{\varphi, Three_s}, t_{\varphi, Five_a}).\]

The interpretations of the types from the extended context are defined by our procedure as
\[\|\T_{\varphi, Three_{s:S}}\| \in \|Three_{s: S}\|\]
\[\|\T_{\varphi, Five_{a:A}}\| \in \|Five_{a: A}\|\]
and moreover
\[\| \T_{\varphi, \frac{Three_{s: S}}{Five_{a: A}}} \| = \|\T_{\varphi, Three_{s:S}}\| \times \ \|\T_{\varphi, Five_{a:A}}\|.\]
\begin{center} \xext=2000 \yext=1300
\begin{picture}(\xext,\yext)(\xoff,\yoff)
\put(-200,100){\framebox(2200,900){}}

  \put(0,300){\framebox(600,500){}}
  \put(0,300){\makebox(600,100){\mbox{$_{\lk John,article 1\rk}$}}}
  \put(0,400){\makebox(600,100){\mbox{$_{\lk John,article 2\rk}$}}}
  \put(0,500){\makebox(600,100){\mbox{$_{\lk John,article 3\rk}$}}}
  \put(0,600){\makebox(600,100){\mbox{$_{\lk John,article 4\rk}$}}}
  \put(0,700){\makebox(600,100){\mbox{$_{\lk John,article 5\rk}$}}}

  \put(600,300){\framebox(600,500){}}
  \put(600,300){\makebox(600,100){\mbox{$_{\lk Bob, article 1\rk}$}}}
  \put(600,400){\makebox(600,100){\mbox{$_{\lk Bob, article 2\rk}$}}}
  \put(600,500){\makebox(600,100){\mbox{$_{\lk Bob, article 3\rk}$}}}
  \put(600,600){\makebox(600,100){\mbox{$_{\lk Bob, article 4\rk}$}}}
  \put(600,700){\makebox(600,100){\mbox{$_{\lk Bob, article 5\rk}$}}}

  \put(1200,300){\framebox(600,500){}}
  \put(1200,300){\makebox(600,100){\mbox{$_{\lk Ken, article 1\rk}$}}}
  \put(1200,400){\makebox(600,100){\mbox{$_{\lk Ken, article 2\rk}$}}}
  \put(1200,500){\makebox(600,100){\mbox{$_{\lk Ken, article 3\rk}$}}}
  \put(1200,600){\makebox(600,100){\mbox{$_{\lk Ken, article 4\rk}$}}}
  \put(1200,700){\makebox(600,100){\mbox{$_{\lk Ken, article 5\rk}$}}}

\put(2100,1150){\makebox(310,100){\mbox{$\|Write\|\subseteq \| S\| \times \|A\|$}}}
\put(2100,1110){\vector(-1,-1){380}}
 \end{picture}
\end{center}
The anaphoric continuation then quantifies universally over the respective interpretations, yielding the desired truth-conditions \textit{All of the three scientists cooperated in presenting at major conferences all of the five articles that they co-authored}

\subsection{`Donkey anaphora'}

Our treatment of `donkey anaphora' does not run into the `proportion problem' and accommodates ambiguities claimed for `donkey sentences'. Consider an example in (4)
\begin{enumerate}[(4)]
\item \label{example4} Every farmer who owns a donkey beats it.
\end{enumerate}
On our analysis, pronouns in `donkey sentences' quantify over (possibly dependent) types introduced by modified common nouns ($MCN$).

To account for the dynamic contribution of modified common nouns, we include in our system {\em $\ast$-sentences} (i.e.\ sentences with dummy quantifier phrases):

\[ \Gamma_{\varphi_1}\vdash \varphi_0 : \textit{Every farmer beats it.} \hskip 12cm\]
\[ \Gamma\vdash \varphi_1 : \textit{Farmer owns a donkey} : \textit{$\ast$-sentence} \hskip 12cm\]

\vspace{2mm}

\noindent The $MCN$ (= \textit{farmer who owns a donkey}) translates into a $\ast$-sentence (with a dummy-quantifier phrase $f:F$)
\[\Gamma \vdash f:F|\exists_{d:D}Own(f, d),\]
and extends the context by adding new variable specifications on newly formed types for every (dummy-) quantifier phrase in the {\em pointed chain} $Ch^{*}$ (= $f:F|\exists_{d:D}$)
\[t_{\varphi, f}: \T_{\varphi,{f:F}}; \  t_{\varphi, \exists_d} : \T_{\varphi, \exists_{d:D}}(t_{\varphi, f}).\]

The interpretations of the types from the extended context $\Gamma_\varphi$ are defined in our usual two-step algorithm. Thus the $\ast$-sentence extends the context by adding new variable specifications on newly formed types
\begin{itemize}
    \item the type $\T_{\varphi,{f:F}}$ interpreted as $\|\T_{\varphi,{f:F}}\|$ (i.e.\ the set of farmers who own some donkeys);
    \item the dependent type $\T_{\varphi, \exists_{d:D}}(t_{\varphi, f})$, interpreted for $a\in \|\T_{\varphi,{f:F}}\|$ as $\|\T_{\varphi, \exists_{d:D}}\|(a)$ (i.e.\ the set of donkeys owned by the farmer $a$).
\end{itemize}
The main clause  $\varphi_0$ (= \textit{Every farmer beats it}) quantifies universally over the respective interpretations:
\[\Gamma_{\varphi_1} \vdash \forall_{t_{\varphi, f}: \T_{\varphi,{f:F}}} | \forall_{t_{\varphi, \exists_d} : \T_{\varphi, \exists_{d:D}}(t_{\varphi, f})} Beat (t_{\varphi, f}, t_{\varphi, \exists_d}),\]
giving the correct truth conditions \textit{Every farmer who owns a donkey beats every donkey he owns.}

Our analysis can be extended to account for more complicated `donkey sentences' such as \textit{Every farmer who owns donkeys beats most of them}. Importantly, the solution does not run into the `proportion problem'. Since we quantify over farmers and the respective fibers of the donkeys owned (and not over $\langle farmer, donkey \rangle$ pairs), a sentence like \textit{Most farmers who own a donkey beat it} comes out false if there are ten farmers who own one donkey and never beat them, and one farmer who owns twenty donkeys and beats all of them. Furthermore, sentences like (4) have been claimed to be ambiguous between the so-called (i) strong reading: \textit{Every farmer who owns a donkey beats \textsc{every} donkey he owns}, and (ii) weak reading: \textit{Every farmer who owns a donkey beats \textsc{at least one} donkey he owns}. Our analysis can accommodate this observation by taking the weak reading to simply employ the quantifier \textit{some} in place of \textit{every} (e.g. we can assume that pragmatic factors (world knowledge, discourse context) can sometimes override the maximality constraint associated with anaphoric pronouns, i.e. under special circumstances, anaphoric pronouns can be treated as existential quantifiers).

\subsection{Nested dependencies}\label{5}

As the type dependencies can be nested, our analysis can be extended to sentences involving three and more quantifiers. Consider an example in (5)
\begin{enumerate}[(5)]
\item\label{example5} Every student bought most professors a flower. They will give them to them tomorrow.
\end{enumerate}
The first sentence in (5)(on the interpretation where \textit{a flower} depends on \textit{most professors} that depends on \textit{every student}) translates into
\[\Gamma \vdash \forall_{s:S} | Most_{p:P}| \exists_{f:F} Buy(s,p,f),\]
and extends the context by adding new variable specifications on newly formed types for every quantifier phrase in $Ch$
\[t_{\varphi, \forall_s}: \T_{\varphi, \forall_{s:S}}; \  t_{\varphi, Most_p} : \T_{\varphi, Most_{p:P}}(t_{\varphi, \forall_s}); \  t_{\varphi, \exists_f} : \T_{\varphi, \exists_{f:F}}(t_{\varphi, \forall_s},t_{\varphi, Most_p}) \]

We now apply our interpretation algorithm.

\noindent \textbf{Step 1.} We define fibers of new types by inverse induction.

\noindent Basic step. For the whole chain $Ch = \forall_{s:S} | Most_{p:P}| \exists_{f:F}$ we put

\[ \|\T_{\varphi, \forall_{s:S} | Most_{p:P}| \exists_{f:F}}\|:= \|Buy\|. \]

\noindent Inductive step.

\small {\[\|\T_{\varphi, \forall_{s:S}}\| = \hskip 11 cm \]
\[= \{a \in \|S\| \, : \,  \{b \in \|P\| \, : \,  \{c \in \|F\| \, : \, \lk a,b,c \rk\in \|Buy\|\} \in  \|\exists_{f:F}\| \} \in  \| Most_{p:P}\| \} \]}

\noindent and for $a\in \|M\|$,

\[\|\T_{\varphi, Most_{p:P}}\|(a) = \{b \in \|P\| \, : \,  \{c \in \|F\| \, : \, \lk a,b,c \rk\in \|Buy\|\} \in  \|\exists_{f:F}\| \} \hskip 8cm\]

\noindent and for $a\in \|M\|$ and $b\in \|P\|$,

\[\|\T_{\varphi, \exists_{f:F}}\|(a,b) = \{c \in \|F\| \, : \, \lk a,b,c \rk\in \|Buy\| \} \hskip 8cm\]

\noindent \textbf{Step 2.} We build dependent types from fibers.

\small {\[\|\T_{\varphi, \forall_{s:S}}\| = \hskip 11 cm \]
\[= \{a \in \|S\| \, : \,  \{b \in \|P\| \, : \,  \{c \in \|F\| \, : \, \lk a,b,c \rk\in \|Buy\|\} \in  \|\exists_{f:F}\| \} \in  \| Most_{p:P}\| \} \]}
\[ \|\T_{\varphi, Most_{p:P}}\| = \bigcup\{ \{ a \}\times \|\T_{\varphi, Most_{p:P}}\|(a) : a \in \|\T_{\varphi, \forall_{s:S}}\|\} \hskip 8cm\]
\[ \|\T_{\varphi, \exists_{f:F}}\| = \bigcup\{ \{ \lk a , b\rk \}\times \|\T_{\varphi, \exists_{f:F}}\|(a,b) : a \in \|\T_{\varphi, \forall_{s:S}}\|, b \in \|\T_{\varphi, Most_{p:P}}\|(a) \} \hskip 8cm\]

Thus the first sentence in (5) extends the context by adding new variable specifications on newly formed types
\begin{itemize}
  \item the type $\T_{\varphi, \forall_{s:S}}$ interpreted as $\|\T_{\varphi, \forall_{s:S}}\|$ (i.e.\ the set of students who bought for most of their professors a flower);
  \item the dependent type $\T_{\varphi, Most_{p:P}}(t_{\varphi, \forall_s})$, interpreted for $a\in \|\T_{\varphi, \forall_{s:S}}\|$ as \\ $\|\T_{\varphi, Most_{p:P}}\|(a)$ (i.e.\ the set of professors for whom the student $a$ bought flowers);
  \item another dependent type $\T_{\varphi, \exists_{f:F}}(t_{\varphi, \forall_s},t_{\varphi, Most_p}) $, interpreted for $a\in \|\T_{\varphi, \forall_{s:S}}\|$ and $b\in \|\T_{\varphi, Most_{p:P}}\|(a)$ as $\|\T_{\varphi, \exists_{f:F}}\|(a,b)$ (i.e.\ the set of flowers that the student $a$ bought for the professor $b$).
\end{itemize}
In the second sentence of (5) the three pronouns $they_s$, $them_p$, and $them_f$ quantify universally over the respective interpretations. The anaphoric continuation in (5) translates into
\[\Gamma_\varphi \vdash \forall_{t_{\varphi, \forall_s}: \T_{\varphi, \forall_{s:S}}}| \forall_{t_{\varphi, Most_p} : \T_{\varphi, Most_{p:P}}(t_{\varphi, \forall_s})}| \forall_{t_{\varphi, \exists_f} : \T_{\varphi, \exists_{f:F}}(t_{\varphi, \forall_s}, t_{\varphi, Most_p})} \] \[Give(t_{\varphi, \forall_s}, t_{\varphi, Most_p}, t_{\varphi, \exists_f}),\]
yielding the correct truth conditions \textit{Every student will give the respective professors the respective flowers he bought for them}.

\subsection{Escaping dependencies}

Unbound anaphoric pronouns are interpreted with reference to the context created by the foregoing text, i.e. they can refer to what is given in the context (referents, dependencies). There are cases, however, where we want pronouns to escape certain dependencies (see \cite{NR}). This is necessary to get the proper reading of the second sentence in (6)

\begin{enumerate}[(6)]
\item Every man loves a woman. They (the women) are (all) smart.
\end{enumerate}
The way to understand the anaphoric continuation is that all of the women loved are smart. The pronoun \textit{they} in the anaphoric continuation refers to the entire set of women loved (by particular men).

On our analysis, the first sentence extends the context by adding new variable specifications on newly formed types
\begin{itemize}
    \item the type $\T_{\varphi, \forall_{m:M}}$, interpreted as $\|\T_{\varphi, \forall_{m:M}}\|$ (i.e.\ the set of men who love some women);
    \item the dependent type $\T_{\varphi, \exists_{w:W}}(t_{\varphi, m})$, interpreted for $a\in \|\T_{\varphi, \forall_{m:M}}\|$ as \\$\|\T_{\varphi, \exists_{w:W}}\|(a)$ (i.e.\ the set of women loved by the man $a$).
\end{itemize}
The pronoun \textit{they} in the anaphoric continuation quantifies universally over the set of ALL women loved, escaping a dependency on the $man$-variable. In our system, this process is enabled by a type constructor $\Sigma$.
\[\Sigma_{t_{\varphi, \forall_m} : \T_{\varphi, \forall_{m:M}}} \T_{\varphi, \exists_{w:W}}(t_{\varphi, \forall_m}), \]
for short
\[\Sigma_{t_m:\T_M} \T_W (t_m) \]
interpreted as
\[ \|\Sigma_{t_m:\T_M} \T_W (t_m)\| \;= \;
                \coprod_{a\in\|\T_M\|} (\{ a \} \times \|\pi_{\T_W,t_m}\|^{-1}(a))\]
i.e. we take the sum of fibers of women over men in $\|\T_M\|.$\\

\noindent The pronoun \textit{they} in the anaphoric continuation quantifies universally over the set $\|\Sigma_{t_m:\T_M(w)} \T_W (t_m)\|$, yielding the correct truth conditions \textit{Every woman loved is smart}.

Consider now a more complicated example (a variant of the example introduced in Section \ref{5})
 \begin{enumerate}[(7)]
\item\label{example7} Every student bought most professors a flower. They picked them carefully.
\end{enumerate}
To get the proper reading of the second sentence: Each student picked carefully all of the flowers bought for most of his professors, we need the second pronoun to escape a dependency on the $professor$-variable. The pronoun $them$ in the anaphoric continuation quantifies universally over the set of ALL flowers that the student $a\in \|\T_{\varphi, \forall_{s:S}}\|$ bought for the professors in $\|\T_{\varphi, Most_{p:P}}\|(a)$

\[\Sigma_{t_{\varphi, Most_p} : \T_{\varphi, Most_{p:P}}(t_{\varphi, \forall_s})} \T_{\varphi, \exists_{f:F}}(t_{\varphi, \forall_s}, t_{\varphi, Most_p}), \]
for short
\[\Sigma_{t_p:\T_P(t_s)} \T_F (t_s, t_p) \]
interpreted as
\[ \|\Sigma_{t_p:\T_P(t_s)} \T_F (t_s, t_p)\|(a) \;= \;
                \coprod_{b\in\|\T_P\|(a)} (\{ b \} \times \|\pi_{\T_F,t_p}\|^{-1}(b))\]
i.e. we take the sum of fibers of flowers over professors for whom the student $a\in \|\T_{\varphi, \forall_{s:S}}\|$  bought flowers.

Thus, in our example, the context gets updated by adding a new variable specification on a newly formed $\Sigma$-type (abbrev. $\T_{\varphi, \Sigma}$)
\[t_{\varphi, \forall_s}: \T_{\varphi, \forall_{s:S}}; \  t_{\varphi, \Sigma}: \T_{\varphi, \Sigma}(t_{\varphi, \forall_s})\]
The anaphoric continuation in (7) translates into
\[\Gamma_\varphi \vdash \forall_{t_{\varphi, \forall_s}: \T_{\varphi, \forall_{s:S}}} | \forall_{t_{\varphi, \Sigma}: \T_{\varphi, \Sigma}(t_{\varphi, \forall_s})} Pick (t_{\varphi, \forall_s}, t_{\varphi, \Sigma}), \]
yielding the correct truth conditions \textit{Every student picked all flowers he bought for most his professors carefully}.

To accommodate all of such extra processes needed to obtain a new context out of the old one we introduce a {\em refresh operation}. The refresh operation will include: addition of variable declarations on presupposed types (where by presupposed types we understand types belonging to the relevant common ground shared by the speaker and hearer); $\sum$, $\prod$ of the types given in the context, etc.

\subsection{Iterated `donkey examples'}
Finally, we will show how our system handles iterated `donkey sentences'. Consider an example in (8)
\begin{enumerate}[(8)]
\item \label{example8} Every hunter who owns a dog who chases a fox helps him get it.
\end{enumerate}
The sentence in (8) quantifies over (possibly dependent) types determined by the {\em type specifying sequence} of $\ast$-sentences, $\vec{\varphi}$
\[  \Gamma_{\vec\varphi}\vdash \varphi_0 : \textit{Every hunter helps him get it.}\hskip 12cm\]
where $\vec{\varphi}=\lk \varphi_1, \varphi_2 \rk$ (linked via the $dog$-variable)
\[  \Gamma\vdash \varphi_1 : \textit{Hunter owns a dog} : \textit{$\ast$-sentence} \hskip 12cm\]
\[  \Gamma\vdash \varphi_2 : \textit{Dog chases a fox} : \textit{$\ast$-sentence} \hskip 12cm\]

Using $\T$-constructor we define the context $\Gamma_{\vec{\varphi}}$ as
\[ \Gamma_{\vec{\varphi}}= \Gamma, \T(\varphi)\]
where $\varphi$ is the $\ast$-sentence
\[  h:H|\exists_{d:D}|\exists_{f:F}\ Own(h,d) \wedge Chase (d,f) \]

The interpretation of the types from the extended context $\Gamma_{\vec{\varphi}}$ are defined in our two-step algorithm.\\
\textbf{Step 1.}\\
\noindent Basic step. For the whole chain we put
\[ \|\T_{h:H|\exists_{d:D}|\exists_{f:F}}\|:= \|Own(h,d)\;\ \wedge \;\ Chase (d,f)\|. \]
Inductive step.
\[\|\T_{\varphi, h:H}\| = \{a \in \|H\| : \{b \in \|D\| :  \{c \in \|F\|  :\]
\[\lk a,b,c \rk\in \|Own(h,d)\;\ \wedge \;\ Chase (d,f)\|\}  \in  \|\exists_{f:F}\| \} \in  \| \exists_{d:D}\| \} \]
and for $a\in \|H\|$,
\[\|\T_{\varphi, \exists_{d:D}}\|(a) = \{b \in \|D\| \, : \,  \{c \in \|F\|  : \]
 \[ \lk a,b,c \rk\in \|Own(h,d)\;\ \wedge \;\ Chase (d,f)\|\} \in  \|\exists_{f:F}\| \} \]
and for $a\in \|H\|$ and $b\in \|D\|$,
\[\|\T_{\varphi, \exists_{f:F}}\|(a,b) = \{c \in \|F\| \, : \, \lk a,b,c \rk\in \|Own(h,d)\;\ \wedge \;\ Chase (d,f)\| \}, \]
where
\[ \|Own(h, d) \;\ \wedge \;\ Chase(d, f)\| = \pi_{h,d}^{-1}(\|Own(h, d)\|) \;\ \cap \;\ \pi_{d,f}^{-1}(\|Chase(d, f)\|), \hskip 8cm\]
with $\pi_{h,d}: \|h\| \times \|d\| \times \|f\| \rightarrow  \|h\| \times \|d\|$ and $\pi_{d,f}: \|h\| \times \|d\| \times \|f\| \rightarrow  \|d\| \times \|f\|.$

\vspace{2mm}
\normalsize
\noindent\textbf{Step 2.} 

\[ \|\T_{\varphi, \exists_{d:D}}\| = \bigcup\{ \{ a \}\times \|\T_{\varphi, \exists_{d:D}}\|(a) : a \in \|\T_{\varphi, h:H}\|\} \hskip 8cm\]
\[ \|\T_{\varphi, \exists_{f:F}}\| = \bigcup\{ \{ \lk a , b\rk \}\times \|\T_{\varphi, \exists_{f:F}}\|(a,b) : a \in \|\T_{\varphi, h:H}\|, b \in \|\T_{\varphi, \exists_{d:D}}\|(a) \} \hskip 8cm\]

\vspace{2mm}

\noindent The main clause  $\varphi_0$ (= \textit{Every hunter helps him get it}) quantifies universally over the respective interpretations,
giving the correct truth conditions.

\section{System - syntax}

This and the following section define, respectively, the syntax and the semantics of our system. As types can depend on variables in our system, we have three kinds of occurrences of variables: binding (next to quantifiers), indexing (next to types), and argument (in quantifier-free formulas). As a variable in a formula might appear in any of those roles, this has to be taken into account when building formulas and defining their semantics. This is the source of the main technical difficulty (the provisos included in definitions) and the increased complication of the system.

\subsection{Alphabet}
The alphabet consists of
\begin{enumerate}
  \item type variables $X, Y, Z, \ldots$;
  \item type constants $M, men, women, \ldots$;
  \item type constructors: $\sum, \prod, \T$;
  \item individual variables $x, y, z, \ldots$;
  \item predicates $P, P',P_1, \ldots$ (with arities specified);
  \item connectives $\wedge$;
  \item quantifier symbols $\exists, \forall, Three, Five, Q_1, Q_2, \ldots$;
  \item three chain constructors: $?|?$, $\frac{\hskip 2mm ?\hskip 2mm}{?}$, $(?,\ldots,?)$.
\end{enumerate}

\subsection{Contexts}

A context is a list of type specifications of (individual) variables. Empty context $\emptyset$ is a context. If we have a context
\[ \Gamma = x_1 : X_1, \ldots, x_k : X_k(\lk x_i\rk_{i\in J_k}), \ldots, x_n : X_n(\lk x_i\rk_{i\in J_n}) \]
then the judgement
\[ \vdash \Gamma : \textsf{context}\]
expresses this fact. Having a context $\Gamma$ as above, we can declare a type $X_{n+1}$ in that context
\[ \Gamma \vdash X_{n+1}(\lk x_i\rk_{i\in J_{n+1}}) : \textsf{type} \]
where $J_{n+1}\subseteq \{ 1, \ldots, n\}$ such that if $i\in J_{n+1}$, then $J_i\subseteq J_{n+1}$, $J_1=\emptyset$.  The type $X_{n+1}$ depends on variables $\lk x_i\rk_{i\in J_{n+1}}$. Now, we can declare a new variable of the type $X_{n+1}(\lk x_i\rk_{i\in J_{n+1}})$ in the context $\Gamma$
\[ \Gamma  \vdash x_{n+1} : X_{n+1}(\lk x_i\rk_{i\in J_{n+1}})\]
and extend the context $\Gamma$ by adding this variable declaration, i.e. we have
\[ \vdash \Gamma,  x_{n+1} : X_{n+1}(\lk x_i\rk_{i\in J_{n+1}}) : \textsf{context}\]
$\Gamma'$ is a {\em subcontext} of $\Gamma$ if $\Gamma'$ is a context and a sublist of $\Gamma$. Let $\Delta$ be a list of variable declarations from a context $\Gamma$, $\Delta'$ the least subcontext of $\Gamma$ containing $\Delta$. We say that $\Delta$ is {\em convex} iff $\Delta'-\Delta$ is again a context.

The variables the types depend on are always explicitly written down in declarations. We can think of a context as (a linearization of) a partially ordered set of declarations such that the declaration of a variable $x$ (of type $X$) precedes the declaration of the variable $y$ (of type $Y$) iff the type $Y$ depends on the variable $x$.

\subsection{Type formation: $\Sigma$-types and $\Pi$-types} \label{2}
Having a type declaration
\[ \Gamma, y : Y(\vec{x})  \vdash  Z (\vec{y})  : \textsf{type}\]
with $y$ occurring in the list $\vec{y}$
we can declare $\Sigma$-type
\[ \Gamma \vdash \Sigma_{y:Y(\vec{x})} Z (\vec{y})  : \textsf{type}\]
and also $\Pi$-type
\[ \Gamma \vdash \Pi_{y:Y(\vec{x})} Z (\vec{y})  : \textsf{type}\]
So declared types do not depend on the variable $y$. Now we can declare new variables of those types. 

\subsection{Quantifier-free formulas} For our purpose we need only atomic formulas and their conjunctions. We have
\[ \Gamma \vdash P(x_1,\ldots,x_n) : \textsf{qf-formula}   \]
whenever $P$ is an $n$-ary predicate and the declarations of the variables $x_1,\ldots,x_n$ form a subcontext of $\Gamma$.
Moreover, we have a formation rule for the conjunction of quantifier-free formulas
\[ \frac{\Gamma \vdash A^i(x^i_1,\ldots,x^i_{n_i}) : \textsf{qf-formula}\;\;{\rm for}\;\; i=1,\ldots, m}
{\Gamma \vdash \bigwedge_{i=1}^mA^i(x^i_1,\ldots,x^i_{n_i}) : \textsf{qf-formula}}\]

\subsection{Quantifier phrases}
If we have a context $\Gamma, y:Y(\vec{x}),\Delta$ and quantifier symbol $Q$, then we can form a {\em quantifier phrase} $Q_{y:Y(\vec{x})}$ in that context. We write
\[ \Gamma, y:Y(\vec{x}),\Delta \vdash Q_{y:Y(\vec{x})} : \textsf{QP}  \]
to express this fact. In a quantifier prase $Q_{y:Y(\vec{x})}$
\begin{enumerate}
  \item the variable $y$ is the {\em binding variable} and
  \item the variables $\vec{x}$ are {\em indexing variables}.
\end{enumerate}

\subsection{Packs of quantifiers}
Quantifiers phrases can be grouped together to form a pack of quantifiers. The pack of quantifiers formation rule is as follows.
  \[ \frac{\Gamma \vdash Q_{i\; y_i:Y_i(\vec{x}_i)} : \textsf{QP} \;\;\; i=1,\ldots k}{\Gamma \vdash (Q_{1\; y_1:Y_1(\vec{x}_1)},\ldots,Q_{k\; y_k:Y_k(\vec{x}_k)}  ): \textsf{pack} } \]
  where, with $\vec{y}=y_1,\ldots, y_k$ and $\vec{x}= \bigcup_{i=1}^k \vec{x}_i$, we have that $y_i\neq y_j$ for $i\neq j$ and $\vec{y}\cap \vec{x}=\emptyset$. In so constructed pack
  \begin{enumerate}
    \item the binding variables are $\vec{y}$ and
    \item the indexing variables are $\vec{x}$.
  \end{enumerate}
  We can denote such a pack $Pc_{\vec{y}:\vec{Y}(\vec{x})}$ to indicate the variables involved. One-element pack will be denoted and treated as  a quantifier phrase. This is why we denote such a pack as $Q_{y:Y(\vec{x})}$ rather than $(Q_{y:Y(\vec{x})})$.

\subsection{Pre-chains and chains of quantifiers} Chains and pre-chains of quantifiers have binding variables and indexing variables. By $Ch_{\vec{y}:\vec{Y}(\vec{x})}$ we denote a pre-chain with binding variables $\vec{y}$ and indexing variables $\vec{x}$ so that the type of the variable $y_i$ is $Y_i(\vec{x}_i)$ with $\bigcup_i\vec{x}_i = \vec{x}$.  Chains of quantifiers are pre-chains in which all indexing variables are bound. Pre-chains of quantifiers arrange quantifier phrases into $N$-free pre-orders, subject to some binding conditions. Mutually comparable QPs in a pre-chain sit in one pack. Thus the pre-chains are built from packs via two chain-constructors of sequential $?|?$ and parallel composition $\frac{?}{?}$.

The chain formation rules are as follows.
\begin{enumerate}
  \item Packs of quantifiers are pre-chains of quantifiers  with the same binding variable and the same indexing variables, i.e.
  \[ \frac{\Gamma \vdash Pc_{\vec{y}:\vec{Y}(\vec{x})} : \textsf{pack}}{\Gamma \vdash Pc_{\vec{y}:\vec{Y}(\vec{x})}: \textsf{pre-chain} } \]
  \item {\em Sequential composition of pre-chains}
  \[  \frac{\Gamma \vdash Ch_{1\; \vec{y}_1:\vec{Y}_1(\vec{x}_1)}: \textsf{pre-chain},\hskip 5mm \Gamma \vdash Ch_{2\; \vec{y}_2:\vec{Y}_2(\vec{x}_2)}: \textsf{pre-chain}}{\Gamma \vdash Ch_{1\; \vec{y}_1:\vec{Y}_1(\vec{x}_1)}|Ch_{2\;\vec{y}_2:\vec{Y}_2(\vec{x}_2)} : \textsf{pre-chain}  }\]
provided
\begin{enumerate}
  \item $\vec{y}_2\cap (\vec{y}_1\cup \vec{x}_1)=\emptyset$,
  \item the declarations of the variables $(\vec{x}_1\cup \vec{x}_2) -(\vec{y}_1\cup \vec{y}_2)$ form a context, a subcontext of $\Gamma$.
\end{enumerate}

In so obtained pre-chain
  \begin{enumerate}
    \item the binding variables are $\vec{y}_1\cup\vec{y}_2$ and
    \item the indexing variables are $\vec{x}_1\cup\vec{x}_2$.
  \end{enumerate}
  \item {\em Parallel composition of pre-chains}
  \[  \frac{\Gamma \vdash Ch_{1\;\vec{y}_1:\vec{Y}_1(\vec{x}_1)}: \textsf{pre-chain},\hskip 5mm \Gamma \vdash Ch_{2\; \vec{y}_2:\vec{Y}_2(\vec{x}_2)}: \textsf{pre-chain}}{\Gamma \vdash \frac{Ch_{1\; \vec{y}_1:\vec{Y}_1(\vec{x}_1)}}{Ch_{2\; \vec{y}_2:\vec{Y}_2(\vec{x}_2)}} : \textsf{pre-chain} }\]
provided $\vec{y}_2\cap (\vec{y}_1\cup \vec{x}_1)=\emptyset=\vec{y}_1\cap (\vec{y}_2\cup \vec{x}_2)$.

As above, in so obtained pre-chain
  \begin{enumerate}
    \item the binding variables are $\vec{y}_1\cup\vec{y}_2$ and
    \item the indexing variables are $\vec{x}_1\cup\vec{x}_2$.
  \end{enumerate}
\end{enumerate}

A pre-chain of quantifiers $Ch_{\vec{y}:\vec{Y}(\vec{x})}$ is a {\em chain} iff $\vec{x}\subseteq \vec{y}$. The following
\[ \Gamma \vdash {Ch}_{\vec{y}:\vec{Y}(\vec{x})}:\textsf{chain} \]
expresses the fact that ${Ch}_{\vec{y}:\vec{Y}(\vec{x})}$ is a chain of quantifiers in the context $\Gamma$.

\subsection{Formulas, sentences and $\ast$-sentences}
The formulas have binding variables, indexing variables and argument variables. We write $\varphi_{\vec{y}:Y(\vec{x})} (\vec{z})$ for a formula with  binding variables $\vec{y}$, indexing variables $\vec{x}$ and argument variables $\vec{z}$. We have the following formation rule for formulas
\[ \frac{\Gamma \vdash A(\vec{z}) : \textsf{qf-formula},\hskip 5mm \Gamma \vdash Ch_{\vec{y}:\vec{Y}(\vec{x})}: \textsf{pre-chain},}{\Gamma \vdash Ch_{\vec{y}:\vec{Y}(\vec{x})}\; A(\vec{z}) : \textsf{formula}}  \]
provided $\vec{y}$ is {\em final} in $\vec{z}$, i.e., $\vec{y}\subseteq \vec{z}$ and the list of variable declarations of $\vec{z}- \vec{y}$ is a subcontext of $\Gamma$. In so constructed formula
  \begin{enumerate}
    \item the binding variables are $\vec{y}$;
    \item the indexing variables are $\vec{x}$;
    \item the argument variables are $\vec{z}$.
  \end{enumerate}

A formula $\varphi_{\vec{y}:Y(\vec{x})}(\vec{z})$ is a {\em sentence } iff $\vec{z}\subseteq \vec{y}$ and $\vec{x}\subseteq \vec{y}$. So a sentence is a formula without free variables, neither argument nor indexing.  The following
\[ \Gamma \vdash \varphi_{\vec{y}:Y(\vec{x})}(\vec{z}):\textsf{sentence} \]
expresses the fact that $\varphi_{\vec{y}:Y(\vec{x})}(\vec{z})$ is a sentence formed in the context $\Gamma$.

We shall also consider some special formulas that we call $\ast$-sentences. A formula $\varphi_{\vec{y}:Y(\vec{x})}(\vec{z})$  is a $\ast$-{\em sentence} if $\vec{x}\subseteq \vec{y}\cup\vec{z}$ but the set $\vec{z}- \vec{y}$ is possibly not empty and moreover the type of each variable in $\vec{z}- \vec{y}$ is {\em constant}, i.e., it does not depend on variables of other types. In such case we consider the set $\vec{z}- \vec{y}$ as a set of biding variables of an additional pack called a {\em dummy pack} that is placed in front of the whole chain $Ch$. The chain `extended' by this dummy pack will be denoted by $Ch^\ast$ and called a {\em pointed chain}. Clearly, if $\vec{z}- \vec{y}$ is empty, there is no dummy pack and the chain $Ch^\ast$ is $Ch$, i.e. sentences are $\ast$-sentences without dummy packs. We write

\[ \Gamma \vdash \varphi_{\vec{y}:Y(\vec{x})}(\vec{z}):\textsf{$\ast$-sentence} \]
to express the fact that $\varphi_{\vec{y}:Y(\vec{x})}(\vec{z})$ is a $\ast$-sentence formed in the context $\Gamma$.

Having formed a $\ast$-sentence $\varphi$, we can form a new context $\Gamma_\varphi$ defined in the section \ref{Gamma-phi}.

\vskip 2mm

\noindent \textsc{\textbf{{Notation}}} For semantics we need some notation for the variables in the $\ast$-sentence. Suppose we have a $\ast$-sentence
\[ \Gamma \vdash Ch_{\vec{y}:Y(\vec{x})}\; A(\vec{z}):\textsf{$\ast$-sentence} \]
We define
\begin{enumerate}
  \item The environment of pre-chain $Ch$: $Env(Ch)=Env(Ch_{\vec{y}:\vec{Y}(\vec{x})})$ - is the context defining variables $\vec{x}-\vec{y}$;
  \item The binding variables of pre-chain $Ch$: $Bv(Ch)=Bv(Ch_{\vec{y}:\vec{Y}(\vec{x})})$ - is the convex set of declarations in $\Gamma$ of the binding variables in $\vec{y}$;
  \item $\env(Ch)=\env(Ch_{\vec{y}:\vec{Y}(\vec{x})})$ - the set of variables in the environment of $Ch$, i.e. $\vec{x}-\vec{y}$;
  \item $\bv(Ch)=\bv(Ch_{\vec{y}:\vec{Y}(\vec{x})})$ - the set of biding variables $\vec{y}$;
  \item The environment of a pre-chain $Ch'$ in a $*$-sentence $\varphi=Ch_{\vec{y}:Y(\vec{x})}\; A(\vec{z})$, denoted $Env_\varphi(Ch')$, is the set of binding variables in all the packs in $Ch^\ast$ that are $<_\varphi$-smaller than all packs in $Ch'$.  Note $Env(Ch')\subseteq Env_\varphi(Ch')$. If $Ch'=Ch_1|Ch_2$ is a sub-pre-chain of the chain $Ch_{\vec{y}:Y(\vec{x})}$, then $Env_\varphi(Ch_2)=Env_\varphi(Ch_1)\cup Bv(Ch_1)$ and $Env_\varphi(Ch_1)=Env_\varphi(Ch')$.
\end{enumerate}

\subsection{Type formation $\T$}\label{Gamma-phi}

Suppose we have constructed a $\ast$-sentence in a context
\[ \Gamma \vdash Ch_{\vec{y}:\vec{Y}(\vec{x})}\; A(\vec{z})  : \textsf{$\ast$-sentence}. \]
In the following we write  $\varphi$ for $Ch_{\vec{y}:\vec{Y}(\vec{x})}\; A(\vec{z})$ for short. We form a context $\Gamma_{\varphi}$ adding one type and one variable declaration for each pack in $Packs_{Ch}$ as follows.

Suppose $\Phi \in Packs_{Ch}$ and $\Gamma'$ is an extension of the context $\Gamma$ such that one variable declaration  $ t_{\Phi',\varphi} : T_{\Phi',\varphi}$ was already added for each (true\footnote{True in the sense that it is not dummy.}) pack  $\Phi'\in Packs_{Ch}$ such that $\Phi'<_{Ch} \Phi$ but not for $\Phi$ yet. Then we declare a new type

$$ \Gamma' \vdash T_{\Phi,\varphi}(\lk t_{\Phi',\varphi}\rk_{\Phi'\in Packs_{Ch},\Phi'<_{Ch}\Phi} ) :\textsf{type}$$
and we extend the context $\Gamma'$ by a declaration of a new variable $t_{\Phi,\varphi}$ of that type
$$ \Gamma', t_{\Phi,\varphi} : T_{\Phi,\varphi}(\lk t_{\Phi',\varphi}\rk_{\Phi'\in Packs_{Ch},\Phi'<_{Ch}\Phi} )  :\textsf{context}$$

The convex set of variable declarations defined this way will be denoted by
\[ \T(\varphi) \]
and the context obtained from $\Gamma$ by adding  the new variables declarations corresponding to all the packs $Packs_{Ch}$ as described above will be denoted by
$$\Gamma_{\varphi}= \Gamma, \T(\varphi).$$
Note that the formula $A(\vec{z})$ need not be a predicate, any other quantifier-free formula with the variables  from $\vec{z}$ would do. We shall use this observation in the next subsection.

\subsection{Type specifying sequences}

A {\em type specifying sequence for a sentence $\psi$ in  a context} $\Gamma$ is a linked sequence of $\ast$-sentences $\vec{\varphi}$  such that
\[  \Gamma_{\vec{\varphi}}\vdash \psi : \textsf{sentence}\]
holds.

A {\em linked sequence of $\ast$-sentences} $\vec{\varphi}=\lk \varphi_1,\ldots, \varphi_n \rk$ is a sequence of $\ast$-sentences such that the last variable in the $i$-th $\ast$-sentence is the variable in the dummy pack of the $i+1$-st $\ast$-sentence, i.e. it is a sequence as displayed below
\[  \Gamma\vdash {z_1:Z_1}| Q^2_{z_2:Z_2}|\ldots| Q^{k_1}_{z_{k_1}:Z_{k_1}}\; A_1(z_{1},\ldots,z_{k_1}) : \textsf{$\ast$-sentence}  \]
\[  \Gamma\vdash {z_{k_1}:Z_{k_1}}| Q^{k_1+1}_{z_{k_1+1}:Z_{k_1+1}}|\ldots| Q^{k_2}_{z_{k_2}:Z_{k_2}}\; A_2(z_{k_1},\ldots,z_{k_2}) : \textsf{$\ast$-sentence}  \]
\[ \ldots \]
\[  \Gamma\vdash {z_{k_{n-1}}:Z_{k_{n-1}}}| Q^{k_{n-1}+1}_{z_{k_{n-1}+1}:Z_{k_{n-1}+1}}|\ldots| Q^{k_n}_{z_{k_n}:Z_{k_n}}\; A_n(z_{k_{n-1}},\ldots,z_{k_n}) : \textsf{$\ast$-sentence}  \]
where $1<k_1<k_2<\ldots< k_n$ and all variables $z_1,\ldots, z_{k_n}$ are different.

Using $\T$-constructor we define the context $\Gamma_{\vec{\varphi}}$ as
\[ \Gamma_{\vec{\varphi}}= \Gamma, \T(\varphi)\]
where $\varphi$ is the $\ast$-sentence
\[  {z_1:Z_1}| Q^2_{z_2:Z_2}|\ldots| Q^{k_n}_{z_{k_n}:Z_{k_n}}\; \bigwedge_{i=1}^n A_i(z_{k_{i-1}},\ldots,z_{k_i}) \]
where $k_0=1$.


We may write
\[ \Gamma\vdash \Gamma_{\vec{\varphi}}\vdash \psi : \textsf{sentence}\]
for short, if $\psi$ is a sentence in context $\Gamma_{\vec{\varphi}}$ that was built from a type specifying sequence $\vec{\varphi}$ for $\psi$.

\subsection{Pure stories}

A {\em pure story} is a sequence of sentences in contexts described by type specifying sequences as follows
\[  \Gamma_1\vdash (\Gamma_1)_{\vec{\varphi}^1;\psi_1}\vdash \psi_1 \]
\[  \Gamma_2=\textsf{refresh}((\Gamma_1)_{\vec{\varphi}^1;\psi_1})\vdash (\Gamma_2)_{\vec{\varphi}^2;\psi_2} \vdash \psi_2 \]
\[ \ldots \]
\[  \Gamma_{n}=\textsf{refresh}((\Gamma_{n-1})_{\vec{\varphi}^{n-1};\psi_{n-1}})\vdash (\Gamma_n)_{\vec{\varphi}^n;\psi_n} \vdash \psi_n \]
so that
\begin{enumerate}
  \item $\Gamma_1$ is the initial context of the pure story;
  \item $(\Gamma_i)_{\vec{\varphi}^i;\psi_i}$ is a context $\Gamma_i$ extended by a type specifying sequence $\vec{\varphi}^i$ for a sentence $\psi_i$; $\psi_i$ is a sentence in the context $(\Gamma_i)_{\vec{\varphi}^i;\psi_i}$;
  \item the context $\textsf{refresh}(((\Gamma_i)_{\vec{\varphi}^i})_{\psi_{i}})$ is a context obtained from $((\Gamma_i)_{\vec{\varphi}^i})_{\psi_{i}}$ by addition of variable declarations on presupposed types (weakening), $\sum$, $\prod$ of these, and  other pragmatic processes to be further studied.
\end{enumerate}

\section{System - semantics}

\subsection{Interpretation of dependent types}
 The context $\Gamma$
\[ \vdash x: X(\ldots), \ldots, y : Y(\ldots, x,\ldots), \ldots, z : Z(\ldots, x, y,\ldots) : \textsf{context}\]
gives rise to a dependence graph. A {\em dependence graph} $DG_\Gamma=(T_\Gamma,E_\Gamma)$ for the context $\Gamma$ has types of $\Gamma$ as vertices and an edge $\pi_{Y,x} : Y \ra X$
for every variable declaration $x:X(\ldots)$ in $\Gamma$ and every type $Y(\ldots, x,\ldots)$ occurring in $\Gamma$ that depends on $x$.

The {\em dependence diagram} for the context $\Gamma$ is an association $\|-\|: DG_\Gamma \ra Set$ to every type $X$ in $T_\Gamma$ a set $\|X\|$ and every edge $\pi_{Y,x} : Y \ra X$ in $E_\Gamma$ a function $\|\pi_{Y,x}\| : \|Y\| \ra \|X\|$, so that whenever we have a triangle of edges in $E_\Gamma$
\begin{center} \xext=700 \yext=650
\begin{picture}(\xext,\yext)(\xoff,\yoff)
  \settriparms[1`1`1;300]
  \putDtriangle(0,0)[Z`Y`X;\pi_{Z,x}`\pi_{Z,y}`\pi_{Y,x}]
 \end{picture}
\end{center}
the corresponding triangle of functions
\begin{center} \xext=700 \yext=650
\begin{picture}(\xext,\yext)(\xoff,\yoff)
  \settriparms[1`1`1;300]
  \putDtriangle(0,0)[\|Z\|`\|Y\|`\|X\|;\|\pi_{Z,x}\|`\|\pi_{Z,y}\|`\|\pi_{Y,x}\|]
 \end{picture}
\end{center}
commutes, i.e.

\[\|\pi_{Z,x}\|=\|\pi_{Y,x}\|\circ\|\pi_{Z,y}\|.\]

 The interpretation of the context $\Gamma$, the {\em parameter space} $\| \Gamma\|$, is the limit\footnote{By this we mean the (categorical) limit of the described (dependence) diagram in the category $Set$ of sets and functions. The notion of a limit used here is the usual category-theoretic notion. In particular, the notion of a parameter space makes sense in any category with finite limits. However, the definition we give in the text is a standard representation of this limit and does not require any knowledge of Category Theory.} of the dependence diagram $\|-\|:DG_\Gamma\ra Set$. More specifically,
\[ \| \Gamma\|= \| x: X(\ldots), \ldots, y : Y(\ldots, x,\ldots), \ldots, z : Z(\ldots, x, y,\ldots) \| = \hskip 45mm\]
\[ = \{ \vec{a} : dom(\vec{a})=\var(\Gamma),\;\;  \vec{a}(z)\in \|Z\|(\vec{a}\lceil \env(Z)),\|\pi_{Z,x}\|(\vec{a}(z))=\vec{a}(x), {\rm for}\]
\[ \hskip 8,5cm   z:Z\;{\rm in}\; \Gamma,\; x\in \env{Z}    \} \]
where $\var(\Gamma)$ denotes variables declared in $\Gamma$ and $\env(Z)$ denotes indexing variables of the type $Z$.

\subsection{Interpretation of $\Sigma$- and $\Pi$-types} \label{3}

 As in this paper we are not going to use $\Pi$-types, we only include the interpretation of a $\Sigma$-type. For
\[ \Gamma \vdash \Sigma_{y:Y(\vec{x})} Z (\vec{y})  : \textsf{type}\]
we define

\[ \|\Sigma_{y:Y(\vec{x})} Z (\vec{y})\| \;= \;
                \coprod_{b\in\|Y\|} (\{ b \} \times \|\pi_{Z,y}\|^{-1}(b)) 
              \]

If a variable $x$ of type $X$ occurs in $\vec{y}$ and $x\neq y$, then we define projection
\[ \| \pi_{\Sigma_{y:Y(\vec{x})} Z (\vec{y}),x} \|:  \|\Sigma_{y:Y(\vec{x})} Z (\vec{y})\|\lra \|X\|   \]
so that
\[  \| \pi_{\Sigma_{y:Y(\vec{x})} Z (\vec{y}),x}\| (b,c) = \|\pi_{Z,x}\|(c) \]
for $b\in \|Y\|$ and $c\in \|\pi_{Z,y}\|^{-1}(b)$.

\subsection{Interpretation of predicates, conjunctions, and quantifier symbols}

Both predicates and quantifiers are interpreted polymorphically.

If we have a predicate $P$ defined in a context $\Gamma$
\[  x_1 : X_1, \ldots, x_n: X_n(\lk x_i\rk_{i\in J_n}) \vdash  P(x_1,\ldots, x_n): \textsf{qf-formula}\]
then, for any interpretation of the context $\|\Gamma\|$, it is interpreted as a subset of its parameter set, i.e.  $\|P\|(\|\Gamma\|)\subseteq \|\Gamma\|$.
The conjunction of the quantifier-free formulas is interpreted as the intersection of the interpretations of these formulas.

{\em Quantifier symbol} $Q$ is interpreted as quantifier $\|Q\|$ i.e. an association to every\footnote{This association can be partial.} set $Z$ a subset $\|Q\|(Z)\subseteq \cP(Z)$.

\subsection{Interpretation of chains of quantifiers}
We interpret QP's, packs, pre-chains, and chains in the environment of a sentence $Env_\varphi$. This is the only case that is needed. We could interpret the aforementioned syntactic objects in their natural environment $Env$ (i.e. independently of any given sentence) but it would unnecessarily complicate some definitions. Thus having a ($\ast$-)sentence $\varphi= Ch_{\vec{y}:Y(\vec{x})}\; A(\vec{z})$ (defined in a context $\Gamma$) and
a sub-pre-chain (QP, pack) $Ch'$, for $\vec{a}\in \|Env_\varphi(Ch')\|$, we define the meaning of
\[ \|Ch'\|(\vec{a}) \]

\noindent \textsc{\textbf{{Notation}}}  Let $\varphi=  Ch_{\vec{y}:\vec{Y}}\; A(\vec{y})$ be a $\ast$-sentence built in a context $\Gamma$, $Ch'$ a pre-chain used in the construction of the ($\ast$)-chain $Ch$. Then $Env_\varphi(Ch')$ is a sub-context of $\Gamma$ disjoint from the convex set $Bv(Ch')$ and  $Env_\varphi(Ch'),Bv(Ch')$  is a sub-context of $\Gamma$. For $\vec{a}\in \| Env_\varphi(Ch')\|$ we define $\| Bv(Ch')\|(\vec{a})$ to be the largest set such that
\[  \{ \vec{a}\} \times \| Bv(Ch')\|(\vec{a}) \subseteq \| Env_\varphi(Ch'),Bv(Ch')\|\]

\subsubsection*{Interpretation of quantifier phrases}

{\em Quantifier phrases.} If we have a quantifier phrase
\[ \Gamma \vdash Q_{y:Y(\vec{x})} : \textsf{QP}\]
and $\vec{a}\in \|Env_\varphi(Q_{y:Y(\vec{x})})\|$, then it is interpreted as $\|Q\|(\|Y\|(\vec{a}))\subseteq \cP(\|Y\|(\vec{a}_{\lceil\vec{x}}))$.

\subsubsection*{Interpretation of packs}
If we have a pack of quantifiers in the sentence $\varphi$
\[  Pc=({Q_1}_{y_1:Y_1(\vec{x}_1)}, \ldots {Q_n}_{y_n:Y_n(\vec{x}_n)})  \]
and $\vec{a}\in \| Env_\varphi(Pc)\|$, then its interpretation with the parameter $\vec{a}$ is
\[ \|Pc\|(\vec{a})=\|({Q_1}_{y_1:Y_1(\vec{x}_1)}, \ldots ,{Q_n}_{y_n:Y_n(\vec{x}_n)})\|(\vec{a}) =\hskip 8cm\]
\[=\{ A \subseteq \prod_{i=1}^{n} \|Y_i\|(\vec{a}\lceil\vec{x}_i)\; :\; \pi_i (A)\in \|Q_i\|(\|Y_i\|(\vec{a}\lceil\vec{x}_i),\; {\rm for}\; i=1,\ldots, n\} \]
where $\pi_i$ is the $i$-th projection from the product.

\subsubsection*{Interpretation of chain constructors}
{\em Parallel composition.} For a pre-chain of quantifiers in the sentence $\varphi$
\[ Ch'=\frac{{Ch_1}_{\vec{y}_1:\vec{Y}_1(\vec{x}_1)}}{{Ch_2}_{\vec{y}_2:\vec{Y}_2(\vec{x}_2)}} \]
and $\vec{a}\in \| Env_\varphi(Ch')\|$, we define
\[ \|\frac{{Ch_1}_{\vec{y}_1:\vec{Y}_1(\vec{x}_1)}}{{Ch_2}_{\vec{y}_2:\vec{Y}_2(\vec{x}_2)}}\|(\vec{a}) = \{ A\times B \; :\; A \in \|{Ch_1}_{\vec{y}_1:\vec{Y}_1(\vec{x}_1)}\|(\vec{a}\lceil\vec{x}_1)\; {\rm and}\;\]
\[\hskip 4,5cm B\in \|{Ch_2}_{\vec{y}_2:\vec{Y}_2(\vec{x}_2)}\|(\vec{a}\lceil\vec{x}_2)\}  \]
\vskip 2mm
\noindent {\em Sequential composition.} For a pre-chain of quantifiers in the sentence $\varphi$
\[ Ch'={Ch_1}_{\vec{y}_1:\vec{Y}_1(\vec{x}_1)}|{Ch_2}_{\vec{y}_2:\vec{Y}_2(\vec{x}_2)} \]
and $\vec{a}\in \| Env_\varphi(Ch')\|$, we define
\[ \|{Ch_1}_{\vec{y}_1:\vec{Y}_1(\vec{x}_1)}|{Ch_2}_{\vec{y}_2:\vec{Y}_2(\vec{x}_2)}\|(\vec{a}) = \{ R\subseteq  \|Bv(Ch')\|(\vec{a}) \; :\; \{\vec{b}\in \|Bv(Ch_1)\|(\vec{a}) \; : \; \hskip 2cm\]
\[ \{\vec{c}\in \|Bv(Ch_2)\|(\vec{a},\vec{b}) \; : \; \lk \vec{b},\vec{c} \rk\in R  \}\in
 \|{Ch_2}_{\vec{y}_2:\vec{Y}_2(\vec{x}_2)}\|{(\vec{a},\vec{b})}\}  \in \|{Ch_1}_{\vec{y}_1:\vec{Y}_1(\vec{x}_1)}\|{(\vec{a})} \}  \]

\subsection{Validity}
A sentence
$$  \Gamma \vdash Ch_{\vec{y}:\vec{Y}}\; A(\vec{y})$$
is true under the above interpretation iff

\[ \|A\|(\|\vec{y}:\vec{Y}\|)\in \|Ch_{\vec{y}:\vec{Y}}\|\]

\subsection{Interpretation of dynamic extensions}

Suppose $\Gamma$ is a context 
and that we have a $\ast$-sentence
\[ \Gamma \vdash Ch_{\vec{y}:\vec{Y}(\vec{x})}\; A(\vec{z}) : \textsf{$\ast$-sentence}.\]
As before, we shall write $\varphi$ for $Ch_{\vec{y}:\vec{Y}(\vec{x})}\; A(\vec{z})$. 
We shall describe the interpretation of the context $\Gamma_{\varphi}$ described in Section \ref{Gamma-phi} that extends the interpretation of the context $\Gamma$.

Thus we are given a dependence diagram $\|-\|^{\Gamma}: DG_{\Gamma} \ra Set$ and we shall define another dependence diagram
\[ \|-\|=\|-\|^{\Gamma_{\varphi}}: DG_{\Gamma_{\varphi}} \ra Set\]
extending $\|-\|^{\Gamma}$ to the context $\Gamma_{\varphi}$. Thus, for $\Phi\in Pack_{Ch^\ast}$ we need to define  $\|\T_{\Phi,\varphi}\|$ and for $\Phi'<_{Ch^\ast} \Phi$ we need to define
\[ \|\pi_{\T_{\Phi,\varphi},t_{\Phi'}}\|: \|\T_{\Phi,\varphi}\|\lra \|\T_{\Phi',\varphi}\|  \]

This will be done in two steps:\\
{\bf Step 1.} (Fibers of new types defined by inverse induction.)\\
We shall define, for the sub-prechains $Ch'$ of $Ch^\ast$ and $\vec{a}\in \|Env_\varphi(Ch')\|$, 
a set
\[ \|\T_{\varphi,Ch'}\|(\vec{a})\subseteq \|Bv(Ch')\|(\vec{a}) \]
This is done using the inductive clauses through which we have defined $Ch^\ast$ but in the reverse direction.

The {\em basic step} is when $Ch'$ is equal to the whole pointed chain $Ch^\ast$. In this case, we put
\[ \|\T_{\varphi,Ch'}\|=\|A\| \]
i.e. we interpret $\T_{\varphi,Ch'}$ as the extension of the whole predicate $A$.

The {\em inductive step}.
Now assume that the set  $\|\T_{\varphi,Ch'}\|(\vec{a})$ is already defined for $\vec{a}\in \| Env_\varphi(Ch')\|$.

\noindent {\em Parallel decomposition.} If we have
\[ Ch'=\frac{{Ch_1}_{\vec{y}_1:\vec{Y}_1(\vec{x}_1)}}{{Ch_2}_{\vec{y}_2:\vec{Y}_2(\vec{x}_2)}} \]
then we define sets
\[ \|\T_{\varphi,Ch_i}\|(\vec{a}\lceil\vec{x}_i) \in \|Ch_i\|(\vec{a}\lceil\vec{x}_i)\]
for $i=1,2$ so that
\[ \|\T_{\varphi,Ch'}\|(\vec{a})=\|\T_{\varphi,Ch_1}\|(\vec{a}\lceil\vec{x}_1,)\times \|\T_{\varphi,Ch_2}\|(\vec{a}\lceil\vec{x}_2) \]
if such sets exist, and these sets $\|\T_{\varphi,Ch_i}\|(\vec{a})$ are undefined\footnote{Such sets might not be determined uniquely if one of them is empty.} otherwise.

\noindent {\em Sequential decomposition.}
If we have
\[ Ch'={Ch_1}_{\vec{y}_1:\vec{Y}_1(\vec{x}_1)}|{Ch_2}_{\vec{y}_2:\vec{Y}_2(\vec{x}_2)} \]
then we put
\[ \|\T_{\varphi,Ch_1}\|(\vec{a}) = \{ \vec{b}\in \|Bv(Ch_1)\|(\vec{a}) : \{ \vec{c} \in \|Bv(Ch_2)\|(\vec{a},\vec{b}):\]
\[\lk \vec{b},\vec{c}\rk \in  \|\T_{\varphi,Ch'}\|(\vec{a})  \} \in \|Ch_2\|(\vec{a},\vec{b})   \} \]
For $\vec{b}\in \|Bv(Ch_1)\|$, we put
\[ \|\T_{\varphi,Ch_2}\|(\vec{a},\vec{b})= \{ \vec{c}\in\|Bv(Ch_2)\|(\vec{a},\vec{b}) :\lk \vec{b},\vec{c}\rk\in \|\T_{\varphi,Ch'}\|(\vec{a},) \} \]
{\bf Step 2.} (Building dependent types from fibers.)

If $\Phi$ is a pack in $Ch^\ast$, $\vec{a}\in \|Env_\varphi(\Phi)\|$, then we put
\[ \|\T_{\varphi,\Phi}\| = \bigcup\{ \{ \vec{a}\}\times \|\T_{\varphi,\Phi}\|(\vec{a}) : \; \vec{a}\in \|Env_\varphi(\Phi)\|,\;\;\hskip 3cm \] \[ \hskip 4,5cm \forall_{\Phi'<_{Ch^\ast} \Phi},\;\; (\vec{a}\lceil\env_\varphi(\Phi'))\in \|\T_{\varphi,\Phi'}\|   \} \]
It remains to define the projections between dependent types.\\ If $\Phi'<_\varphi \Phi$, we define
\[ \pi_{\T_{\varphi,\Phi},t_{\varphi,\Phi'} } : \|\T_{\varphi,\Phi}\| \lra \|\T_{\varphi,\Phi'}\| \]
so that
\[ \vec{a} \mapsto \vec{a}\lceil(\env_\varphi(\Phi')\cup \bv{\Phi'}). \]

\section{Conclusion}

It was our intention in this paper to show that adopting a new type-theoretic approach (with dependent types) to generalized quantification allows a natural and elegant treatment of a wide array of anaphoric data involving natural language quantification. The main technical contribution of our paper consists in combining generalized quantifiers with dependent types. Empirically, our system allows a uniform account of both maximal anaphora to quantifiers and the notoriously difficult cases such as quantificational subordination, cumulative and branching continuations, and `donkey anaphora'.

\section*{Acknowledgments}

The work of Justyna Grudzi\'{n}ska was funded by the National Science Center on the basis of decision DEC-2012/07/B/HS1/00301. The authors would like to thank the anonymous reviewers for valuable comments on an earlier version of this paper. The present paper is a modified and largely extended version of this earlier version which appeared in \textit{Proceedings of EACL 2014 Type Theory and Natural Language Semantics Workshop} (\cite{GZ}).

\end{document}